\documentclass{amsart}
\usepackage{amssymb}

\usepackage{graphicx}

\input{tcilatex}
\input tcilatex

\begin{document}
\title[Lamperti's mbp]{Revisiting John Lamperti's maximal branching process}
\author{Thierry Huillet$^{1}$, Servet Martinez$^{2}$}
\address{$^{1}$Laboratoire de Physique Th\'{e}orique et Mod\'{e}lisation\\
Universit\'{e} de Cergy-Pontoise\\
CNRS UMR-8089\\
Site de Saint Martin\\
2 avenue Adolphe-Chauvin\\
95302 Cergy-Pontoise, France\\
$^{2}$Depto. Ingenieria Matematica and Centro Modelamiento Matematico\\
Universidad de Chile\\
UMI 2807, Uchile-Cnrs\\
Casilla 170-3 Correo 3\\
Santiago, Chile\\
E-mail: huillet@u-cergy.fr, smartine@dim.uchile.cl}
\maketitle

\begin{abstract}
Lamperti's maximal branching process is revisited, with emphasis on the
description of the shape of the invariant measures in both the recurrent and
transient regimes. A truncated version of this chain is exhibited,
preserving the monotonicity of the original Lamperti chain supported by the
integers. The Brown theory of hitting times applies to the latter chain with
finite state-space, including sharp strong time to stationarity. Additional
information on these hitting time problems are drawn from the
quasi-stationary point of view.\newline

\textbf{Running title: }Lamperti's MBP.\newline


\end{abstract}

\section{Introduction}

The Lamperti's maximal branching process (mbp) is a modification of the
Galton-Watson (GW) branching process selecting at each step the descendants
of the most prolific ancestor, \cite{L1}. As a Markov chain on the full set
of non-negative integers, Lamperti (\cite{L1}-\cite{L2}) gave sharp
conditions on the tails of the branching number under which this process is
recurrent (either positive or null) or transient.

Our contribution is to describe the corresponding shape of the invariant
measures and we proceed as follows: while fixing a target invariant measure
(supported by the integers) of the mbp, we show (in Proposition $2$) how to
compute in general the law of the branching mechanism that gives rise to it.
Several classes of distributions are supplied both in the recurrent and
transient setups. In Propositions $3$, $4$ and $5$, the target invariant
measures are probabilities with tails getting larger and larger, ranging
from geometric, power-law with index $\alpha \in \left( 0,1\right) $ and
power-law with index $0$ (the target has no moments of any positive order).
In Propositions $6$ (and $7$), it is shown that the null recurrent
(respectively transient) Lamperti chain has a non trivial invariant infinite
and positive measure.

An important feature of the Lamperti chain we also emphasize on is its
failure rate monotonicity (Proposition $1$).

The Lamperti's mbp also makes sense when the branching mechanism takes
values in the finite subset $\left\{ 1,...,N\right\} $ and the question of
computing the law of the branching mechanism giving rise to any finitely
supported target distribution makes sense. We address this point in
Proposition $8$. If the target distribution is in particular the restriction
to $\left\{ 1,...,N\right\} $ of the invariant measure of a mbp with full
state-space, this construction allows to design a truncated version of the
latter chain preserving its failure rate monotonicity feature (Proposition $%
9 $ and Corollary $10$). For failure rate monotone Markov chains with finite
state-space, Brown, \cite{Brown}, designed a theory of hitting times which
thus applies to the truncated Lamperti chain. The main concern is the
relationship existing between the first hitting times of both state $\left\{
N\right\} $ and the restricted invariant measure of the truncated Lamperti
chain. By monotonicity, state $\left\{ N\right\} $ is the largest possible
value that the truncated chain can explore. Under some technical condition
on the initial distribution, it is recalled that the former hitting time
exceeds stochastically the latter (Proposition $11$) which has the structure
of a compound geometric random variable (Proposition $13$). The excess time
is a sharp strong time to stationarity allowing to estimate the distance
between the current state of the truncated chain to its equilibrium
distribution. Its cumulated probability mass function up to $n$ can be
computed from the probability that the truncated chain is in state $\left\{
N\right\} $ after $n$ steps, (Proposition $12$). The alternative classical
quasi-stationary point of view to this problem is also addressed. In
Proposition $14$, we exhibit the rate of decrease of the hitting times to
state $\left\{ N\right\} $ in terms of the quasi-stationary distribution. In
Proposition $15$, we show that under Brown's conditions on the initial
distribution $\mathbf{\pi }_{0}$, the ratio of the large tail probabilities
for the first hitting times of state $\left\{ N\right\} $ starting from $%
\mathbf{\pi }_{0}$ against the quasi-stationary distribution exceeds $1$.
Proposition $16$ deals with a question raised by Brown concerning asymptotic
exponentiality of the hitting times which applies to the truncated Lamperti
chain and its time-reversal.

\section{Lamperti's model}

The Lamperti maximal branching process (mbp) process may be described as an
extremal analogue of the GW branching process, where the next generation is
formed by the offspring of a most productive individual, \cite{L1}. As a
result of some selection (or detection) mechanism, iteratively in each
generation, only the offspring of one of the most productive individuals of
the underlying GW process with branching number $\nu $ is kept (or
detected), the other ones being wiped out (or missed by the detector). This
output mechanism amounts to pruning Galton-Watson trees by iterative
selection of a largest family size ending up with the sub-tree of the
fittest individuals. In \cite{L1}, Lamperti relates this model to a
percolation problem.

With $X_{n}$ the size of such a population at generation $n$, $F_{n}\left(
j\right) =\mathbf{P}\left( X_{n}\leq j\right) $ and $\nu _{j,n+1}\overset{d}{%
=}\nu $ for all $j$, the dynamics under concern is 
\begin{equation*}
X_{n+1}=\max_{j=1,...,X_{n}}\nu _{j,n+1}\Rightarrow F_{n+1}\left( j\right)
=\sum_{i\geq 0}\mathbf{P}\left( X_{n}=i\right) \mathbf{P}\left( \nu \leq
j\right) ^{i}=\mathbf{E}z^{X_{n}}\mid _{z=\mathbf{P}\left( \nu \leq j\right)
}.
\end{equation*}
with initial condition: $X_{0}\overset{d}{\sim }\mathbf{\pi }_{0}$ with $%
\mathbf{P}\left( X_{0}\leq j\right) :=F_{0}\left( j\right) .$

We denote $\mathbf{E}\left( X_{n+1}\mid X_{n}=i\right) =\mathbf{E}%
\max_{j=1,...,i}\nu _{j}=\mathbf{E}\left( m_{i}\right) $ where $%
m_{i}=\max_{j=1,...,i}\nu _{j}.$

Let $p\left( j\right) :=\mathbf{P}\left( \nu =j\right) $. We will assume
that the set $\left\{ j:p\left( j\right) >0\right\} $ is either $\Bbb{N}%
_{0}:=\left\{ 0,1,2,...\right\} $ or $\Bbb{N}:=\left\{ 1,2,...\right\} $
but, as we shall see, the finite case when $\left\{ j:p\left( j\right)
>0\right\} =\left\{ 1,...,N\right\} $ for some integer $N\gg 1$, will also
be of interest.

We shall let $\phi \left( z\right) =\mathbf{E}z^{\nu }$ be the probability
generating function (pgf) of $\nu .$

We shall distinguish two regimes for the branching number $\nu $:

\subsection{Branching number\textbf{\ }$\nu >0$\textbf{.}}

If $\nu >0$ ($p\left( 0\right) =\mathbf{P}\left( \nu =0\right) =0$ and $%
\mathbf{E}\left( \nu \right) >1$), then $X_{n}>0,$ $\forall n\geq 0$ ($%
X_{0}=1$), owing to 
\begin{equation*}
F_{n+1}\left( 0\right) =\mathbf{P}\left( X_{n+1}=0\right) =\mathbf{E}%
z^{X_{n}}\mid _{z=p\left( 0\right) =0}=\mathbf{P}\left( X_{n}=0\right) =0.
\end{equation*}
We can omit state $0$, being disconnected. One main concern in this context
is whether $X_{n}\rightarrow \infty $ with probability (wp) $1$ (a case of
transience) or to some limiting random variable (rv) $X_{\infty }$ (a case
of recurrence): the tails of $\nu $ matter to decide. In the recurrent case,
what is the shape of the invariant probability measure? In the
null-recurrent and transient cases, what are the shapes of the invariant
measure (no longer probability measures). In particular how are the tails of
the invariant measure related to the tails of $\nu $.\newline

- \textbf{Transition matrix of }$\left\{ X_{n}\right\} $. With $F\left(
j\right) =\mathbf{P}\left( \nu \leq j\right) $, $j\geq 1$, $\left\{
X_{n}\right\} $ is a time-homogeneous Markov chain (MC) on $\Bbb{N}$ with
transition matrix ($\sum_{j\geq 1}P\left( i,j\right) =1-F\left( 0\right)
^{i}=1$) 
\begin{equation*}
P\left( i,j\right) =F\left( j\right) ^{i}-F\left( j-1\right) ^{i}\text{, }%
i,j\geq 1
\end{equation*}
equivalently 
\begin{eqnarray*}
\mathbf{P}\left( X_{n+1}>i\mid X_{n}=i\right) &=&1-F\left( i\right) ^{i} \\
\mathbf{P}_{X_{n}}\left( X_{n+1}>X_{n}\right) &=&1-F\left( X_{n}\right)
^{X_{n}}.
\end{eqnarray*}
Note $P\left( 1,j\right) =\mathbf{P}\left( \nu =j\right) .$\newline

- \textbf{Some properties of} $\left\{ X_{n}\right\} $:

- The Lamperti chain clearly is irreducible and aperiodic.

- It holds that $\mathbf{P}\left( X_{n+1}\leq j\mid X_{n}=i\right)
=:P^{c}\left( i,j\right) =F\left( j\right) ^{i}$ is a decreasing function of 
$i$, for all $j$: the Lamperti MC $\left\{ X_{n}\right\} $ is stochastically
monotone (SM). Equivalently, with $\left\{ >j\right\} $ denoting the upper
set $\left\{ j+1,...\right\} ,$ $\mathbf{P}\left( X_{n+1}>j\mid
X_{n}=i\right) =:P\left( i,\left\{ >j\right\} \right) $ is an increasing
function of $i$, for all $j$ and by induction $P^{n}\left( i,\left\{
>j\right\} \right) $ is an increasing function of $i$, for all $j$ and $n.$
In fact, it has a stronger monotonicity feature:

\begin{proposition}
The Lamperti Markov chain $\left\{ X_{n}\right\} $ is failure-rate monotone.
\end{proposition}

\emph{Proof:} The cumulated transition matrix : $P^{c}\left( i,j\right)
:=\sum_{k=1}^{j}P\left( i,k\right) =:$ $P\left( i,\left\{ \leq j\right\}
\right) $ satisfies:

\begin{equation*}
P^{c}\left( i_{1},j_{1}\right) P^{c}\left( i_{2},j_{2}\right) \geq
P^{c}\left( i_{1},j_{2}\right) P^{c}\left( i_{2},j_{1}\right) ,
\end{equation*}
for all $i_{1}<i_{2}$ and $j_{1}<j_{2}$ (the matrix $P^{c}$ is totally
positive of order $2$, viz TP$_{2}$): the MC $\left\{ X_{n}\right\} $ is
failure rate monotone. Since if $P^{c}$ is TP$_{2}$, $P^{c}\left( i,j\right) 
$ is a decreasing function of $i$, for all $j$, (set $j_{2}=\infty $ in the
last inequality to get $P^{c}\left( i_{1},j_{1}\right) \geq P^{c}\left(
i_{2},j_{1}\right) $), TP$_{2}$ matrices $P^{c}$ form a subclass of SM
matrices $P^{c}$. $\Box $\newline

- \textbf{Generation:} As for all Markov chains, with $\left( \mathcal{U}%
_{n};n\geq 1\right) $ a sequence of independent identically distributed
(iid) uniform-$\left( 0,1\right) $ rvs: 
\begin{equation*}
X_{n+1}=\sum_{j\geq 1}j\cdot \mathbf{1}\left( \mathcal{U}_{n+1}\in \left[
P^{c}\left( X_{n},j-1\right) ,P^{c}\left( X_{n},j\right) \right) \right) .
\end{equation*}
We can also check that, with $F^{-1}\left( y\right) =\inf \left( x:F\left(
x\right) \geq y\right) $ the inverse function of $F$, one has $%
X_{n+1}=F^{-1}\left( \mathcal{U}_{n+1}^{1/X_{n}}\right) .$\newline

- \textbf{Transience versus recurrence:} Note that if $\mathbf{P}\left( \nu
>i\right) \sim \lambda /i$, $\lambda >0$, for large $i$ ($\sim $ meaning
that the ratio of the two terms appearing to the left and right of this
symbol tend to $1$ as $i\rightarrow \infty $), $\mathbf{P}\left(
X_{n+1}>i\mid X_{n}=i\right) \sim 1-e^{-\lambda }>0.$ In this case, 
\begin{equation*}
\mathbf{P}\left( X_{n+1}\leq \left[ ix\right] \mid X_{n}=i\right) \sim
F\left( ix\right) ^{i}\sim e^{-\lambda /x}.
\end{equation*}
and with $Z_{n}=\log X_{n}$%
\begin{equation*}
\mathbf{P}\left( Z_{n+1}-Z_{n}\leq z\mid Z_{n}=\log i\right) \sim
e^{-\lambda e^{-z}},
\end{equation*}
independent of $i$. This shows that for large $i$ and for this choice of $%
\nu $, $\left\{ Z_{n}\right\} $ resembles a random walk with independent
increments whose common law is a Gumbel distribution with mean $m=\log
\lambda +\gamma $ ($\gamma $ the Euler constant). So $\left\{ Z_{n}\right\} $
(and $\left\{ X_{n}\right\} $) drifts to $\infty $ if $\lambda >e^{-\gamma }$
($m>0$) and the basic results of Lamperti in \cite{L1}, \cite{L2} are: 
\begin{equation}
\text{If }\lim \inf_{i}i\mathbf{P}\left( \nu >i\right) <c:=e^{-\gamma }\text{%
, then }X_{n}\overset{a.s.}{\rightarrow }X_{\infty }\text{ (ergodicity),}
\label{L1}
\end{equation}
where $X_{\infty }$ is a non-degenerate rv and ergodicity means positive
recurrence and aperiodicity. 
\begin{equation}
\text{If }\lim \sup_{i}i\mathbf{P}\left( \nu >i\right) >c:=e^{-\gamma }\text{%
, then }X_{n}\rightarrow \infty \text{ wp }1\text{ (transience).}  \label{L2}
\end{equation}
In particular, if $\nu $ has tails heavier than $1/i$ ($i\mathbf{P}\left(
\nu >i\right) \rightarrow \infty $), then $X_{n}\rightarrow \infty $ wp $1,$
(transience). 
\begin{equation}
\begin{array}{l}
\text{Critical case, \cite{L2}: } \\ 
\text{If }\mathbf{P}\left( \nu >i\right) \sim e^{-\gamma }/i+d/\left( i\log
i\right) \text{, the process }\left\{ X_{n}\right\} \text{ is:} \\ 
\text{- positive recurrent if }d<-e^{-\gamma }\pi ^{2}/12 \\ 
\text{- null recurrent if }d\in \left[ -e^{-\gamma }\pi ^{2}/12,e^{-\gamma
}\pi ^{2}/12\right)  \\ 
\text{- transient if }d>e^{-\gamma }\pi ^{2}/12.
\end{array}
\label{L3}
\end{equation}

The case $d=e^{-\gamma }\pi ^{2}/12$ is left open and would require
additional information on the tails of $\nu $ to decide whether here $%
\left\{ X_{n}\right\} $ is transient or null recurrent$.$

Whenever the process $\left\{ X_{n}\right\} $ is ergodic, with $\Phi
_{\infty }\left( z\right) :=\mathbf{E}z^{X_{\infty }}$, the functional
equation 
\begin{equation}
F_{\infty }\left( j\right) =\mathbf{P}\left( X_{\infty }\leq j\right) =\Phi
_{\infty }\left( \mathbf{P}\left( \nu \leq j\right) \right) \text{, }j\geq 1
\label{FE0}
\end{equation}
admits a unique solution for the pair $\left( \mathbf{P}\left( X_{\infty
}\leq j\right) ,\mathbf{P}\left( \nu \leq j\right) \right) $. Because $\Phi
_{\infty }\left( z\right) $ is a pgf with $\Phi _{\infty }\left( 0\right) =0$%
, we have $\Phi _{\infty }\left( z\right) <z$ and so $X_{\infty }$ is
stochastically larger than $\nu $:

\begin{equation}
\text{For all }j\geq 1\text{: }\mathbf{P}\left( X_{\infty }\leq j\right) <%
\mathbf{P}\left( \nu \leq j\right) .  \label{SD}
\end{equation}

Clearly, the maximal branching process asymptotically selects a family size $%
X_{\infty }$ which is larger than the typical family size $\nu $ of the
underlying Galton-Watson process. It is then of utmost interest to solve the
functional equation (\ref{FE0}). As we shall see, the position we will adopt
is the following: suppose one has some initial guess of the limiting rv $%
X_{\infty }$, we will identify the branching number $\nu $ of the Lamperti
mbp realizing this task.

An additional problem of interest: how long does it take for $\left\{
X_{n}\right\} $ to reach $X_{\infty }?$ To have an insight on this question,
we shall ask how long it takes, for a suitably truncated version $%
X_{n}^{\left( N\right) }$ of $X_{n}$, to reach height $N\gg 1$, which is
intuitively more demanding than reaching the invariant measure of the
truncated chain itself. We shall address these points.\newline

- \textbf{Time spent in the worst state. }Whenever the process $\left\{
X_{n}\right\} $ is ergodic, it visits infinitely often all the states, in
particular the state $\left\{ 1\right\} $, and a sample path of it is made
of iid successive non-negative excursions through that state$.$ State $%
\left\{ 1\right\} $ is the worst case of the selection mechanism that the
Lamperti chain realizes. By the ergodic theorem, the fraction of time spent
by $\left\{ X_{n}\right\} $ in this state is $\pi \left( 1\right) =\mathbf{P}%
\left( X_{\infty }=1\right) .$ The expected first return time ($\tau _{1,1}$%
) to state $\left\{ 1\right\} $ is $\mathbf{E}\left( \tau _{1,1}\right)
=1/\pi \left( 1\right) $.

Suppose $\left\{ X_{n}\right\} $ enters state $\left\{ 1\right\} $ from
above at some time $n_{1}.$ The first return time $\tau _{1,1}:=\inf \left(
n>n_{1}:X_{n}=1\mid X_{n_{1}}=1\right) $ to state $\left\{ 1\right\} $ is:

- either $1$ if $X_{n_{1}}$ stays there with probability $P\left( 1,1\right)
=F\left( 1\right) $ in the next step; this corresponds to a trivial
excursion of length $1$ and height $0$.

- or, with probability $1-F\left( 1\right) $, $\left\{ X_{n}\right\} $
starts a true excursion with positive height and length $\tau _{1,1}^{+}\geq
2.$

Thus

\begin{equation*}
\mathbf{E}\left( \tau _{1,1}\right) =\frac{1}{\pi \left( 1\right) }=F\left(
1\right) +\left( 1-F\left( 1\right) \right) \mathbf{E}\left( \tau
_{1,1}^{+}\right) \text{ and}
\end{equation*}
\begin{equation*}
\mathbf{E}\left( \tau _{1,1}^{+}\right) =\frac{1}{1-F\left( 1\right) }\left( 
\frac{1}{\pi \left( 1\right) }-F\left( 1\right) \right) >2,
\end{equation*}
entailing the relationship: $\frac{1}{\pi \left( 1\right) }>2-p\left(
1\right) $. Given $\left\{ X_{n}\right\} $ enters state $\left\{ 1\right\} $
from above at some time $n_{1}$, it stays there with probability $P\left(
1,1\right) =F\left( 1\right) $ in the next step, so $\left\{ X_{n}\right\} $
will quit state $\left\{ 1\right\} $ at time $n_{1}+G$ where $G$ is a
shifted geometric random time with success probability $1-F\left( 1\right) .$
After $n_{1}+G$, the chain moves up before returning to state $\left\{
1\right\} $ again and the time it takes is $\tau _{1,1}^{+}$. Considering
two consecutive instants where $\left\{ X_{n}\right\} $ enters state $%
\left\{ 1\right\} $ from above (defining an alternating renewal process),
the fraction of time spent in state $\left\{ 1\right\} $ is: 
\begin{equation*}
\rho =\frac{\mathbf{E}\left( G\right) }{\mathbf{E}\left( G\right) +\mathbf{E}%
\left( \tau _{1,1}^{+}\right) }.
\end{equation*}
From the expression $\mathbf{E}\left( G\right) =F\left( 1\right) /\left(
1-F\left( 1\right) \right) $ and the value of $\mathbf{E}\left( \tau
_{1,1}^{+}\right) $, we get: 
\begin{equation*}
\rho =F\left( 1\right) \pi \left( 1\right) .\text{ }
\end{equation*}
\newline

- \textbf{Time reversal:}

Suppose $\left\{ X_{n}\right\} $ is ergodic. Let $\pi \left( j\right) =%
\mathbf{P}\left( X_{\infty }=j\right) $, $j\geq 1$. With $\mathbf{\pi }%
^{\prime }=\left( \pi \left( 1\right) ,\pi \left( 2\right) ,...\right) $ the
transpose of the column-vector $\mathbf{\pi }$, $P^{\prime }$ the transpose
of $P$ and $\pi \left( i\right) =\mathbf{P}\left( X_{\infty }=i\right) $ the
stochastic matrix 
\begin{equation*}
\overleftarrow{P}=D_{\mathbf{\pi }}^{-1}P^{\prime }D_{\mathbf{\pi }}
\end{equation*}
is the transition matrix of the time-reversed chain $\left\{
X_{n}^{\leftarrow }\right\} $. Since $\overleftarrow{P}\neq P$, there is no
detailed balance. The process $\left\{ X_{n}^{\leftarrow }\right\} $ is such
that its time-reversal $\left( X_{n}^{\leftarrow }\right) ^{\leftarrow
}=X_{n}$ is stochastically monotone. The backward process $\left\{
X_{n}^{\leftarrow }\right\} $ can be generated as follows, with a
time-reversal flavor: with $\left( J_{n};n\geq 1\right) $ an iid sequence
with $J_{1}\overset{d}{\sim }\mathbf{\pi }$, independent of the $\nu $'s,
consider the Markovian dynamics 
\begin{equation}
Y_{n+1}=J_{n+1}\cdot \mathbf{1}\left( \max_{k=1,...,J_{n+1}}\nu
_{k,n+1}=Y_{n}^{{}}\right) ,  \label{TR}
\end{equation}
giving $Y_{n+1}$ as a $\mathbf{\pi }-$mixture of the number of ancestors
whose most productive individuals produce exactly $Y_{n}^{{}}$ descendants
in a Galton-Watson process with branching number $\nu $. We have

\begin{eqnarray*}
\mathbf{P}\left( Y_{n+1}=j\mid Y_{n}^{{}}=i\right) &=&\pi \left( j\right) 
\mathbf{P}\left( \max_{k=1,...,j}\nu _{k,n+1}=i\right) \\
&=&\pi \left( j\right) \left[ F\left( i\right) ^{j}-F\left( i-1\right)
^{j}\right] =\pi \left( j\right) \mathbf{P}\left( X_{n+1}=i\mid
X_{n}=j\right) ,
\end{eqnarray*}
equivalently 
\begin{equation*}
Q=P^{\prime }D_{\mathbf{\pi }}
\end{equation*}
where $Q\left( i,j\right) =\mathbf{P}\left( Y_{n+1}=j\mid
Y_{n}^{{}}=i\right) .$ The process $Y_{n}^{{}}$ is substochastic (there is a
positive probability that given $Y_{n}^{{}}$ no such index $Y_{n+1}$ exists)
and a coffin state can be added to the state-space $\Bbb{N}$ where the
system is sent to if $Y_{n+1}$ does not exist. Let $\tau _{i}$ be the first
hitting time of the coffin state for $Y_{n}^{{}}$ started at $i$ with $%
\mathbf{P}\left( \tau _{i}=1\right) =1-\pi \left( i\right) ,$ the mass
defect in state $i$ of $Q$. Then $X_{n}^{\leftarrow }=Y_{n}^{{}}\mid \tau
_{Y_{n}^{{}}}>1$ (upon conditioning $Y_{n}^{{}}$ stepwise on the event that
the hitting time of the coffin state exceeds one time unit). The process $%
\left\{ X_{n}^{\leftarrow }\right\} $ thus constructed has the transition
matrix $\overleftarrow{P}$, as required.

\subsection{Branching number\textbf{\ }$\nu \geq 0$\textbf{.}}

If $p\left( 0\right) =\mathbf{P}\left( \nu =0\right) >0$ $:$ the above
functional equation must be considered for $j\geq 0.$

We have $F_{n+1}\left( 0\right) =\sum_{i}\mathbf{P}\left( X_{n}=i\right) 
\mathbf{P}\left( \nu =0\right) ^{i}=\mathbf{E}z^{X_{n}}\mid _{z=p\left(
0\right) }>0.$ At each $n$, there is a positive probability that $X_{n}=0.$
If for some $n$, $X_{n}=0$, clearly $X_{n^{\prime }}=0$ for all $n^{\prime
}>n:$ state $0$ is absorbing. $\left\{ X_{n}\right\} $ is again a Markov
chain now on $\mathbf{N}_{0}$ with transition probability matrix 
\begin{equation}
P\left( i,j\right) =F\left( j\right) ^{i}-F\left( j-1\right) ^{i}\text{, }%
i,j\geq 0  \label{P1}
\end{equation}
in particular with $P\left( i,0\right) =F\left( 0\right) ^{i}>0$.

Two cases arise:

$\left( a\right) $ If $\mathbf{E}\left( \nu \right) \leq 1,$ there is almost
sure (a.s.) extinction of the underlying branching process, say at $\tau _{%
\mathbf{\pi }_{0},0},$ and also therefore of $\left\{ X_{n}\right\} $ at $%
\tau _{\mathbf{\pi }_{0},0}^{X}\leq \tau _{\mathbf{\pi }_{0},0}$. We have $%
\mathbf{P}\left( X_{n}=0\right) =\mathbf{P}\left( \tau _{0}^{X}\leq n\right)
\rightarrow 1$ or $\mathbf{P}\left( X_{n}=0\right) =1$, $\forall n\geq \tau
_{0}^{X}$ ($\tau _{\mathbf{\pi }_{0},0}^{X}$ is the absorption time of $%
\left\{ X_{n}\right\} $ at $0$). In this case, $\Phi _{\infty }\left(
z\right) =1$ for all $z\in \left[ 0,1\right] $ and one possible solution to
the functional equation is $F_{\infty }\left( j\right) =1$, $j\geq 0.$ The
only problem here is to fix the law of $\tau _{\mathbf{\pi }_{0},0}^{X}$
which (with $\mathbf{e}_{0}^{\prime }=\left( 1,0,0,...\right) $ with $1$ in
position $0$), is: 
\begin{equation*}
\mathbf{P}\left( X_{n}=0\right) =\mathbf{P}\left( \tau _{0}^{X}\leq n\right)
=\mathbf{\pi }_{0}^{\prime }P^{n}\mathbf{e}_{0}\text{.}
\end{equation*}

$\left( b\right) $ If $\mathbf{E}\left( \nu \right) >1,$ there is extinction
of the underlying branching process with probability $0<\rho _{e}<1$ ($\rho
_{e}$ the smallest solution in $\left[ 0,1\right] $ to $\phi \left( z\right)
=z$) entailing:

- a.s. extinction: given the underlying branching process certainly goes
extinct (an event with probability $\rho _{e}$), the branching process is
generated by the branching number $\nu _{e}$ with $\mathbf{E}\left( z^{\nu
_{e}}\right) =\phi \left( z\rho _{e}\right) /\rho _{e}$ and $\mathbf{E}%
\left( \nu _{e}\right) \leq 1$, entailing: $X_{n}^{e}\rightarrow 0$ with
probability (wp) $1$. The question is how fast and we are back to the case $%
\left( a\right) $.

- a.s. explosion: given the underlying branching process certainly explodes
(an event wp $1-\rho _{e}$), the branching process is generated by $\nu _{%
\overline{e}}$ characterized by $\mathbf{E}\left( z^{\nu _{\overline{e}%
}}\right) =\left( \phi \left( \rho _{e}+z\left( 1-\rho _{e}\right) \right)
-\rho _{e}\right) /\left( 1-\rho _{e}\right) $ with $\mathbf{P}\left( \nu _{%
\overline{e}}=0\right) =0$ and $\mathbf{E}\left( \nu _{\overline{e}}\right) =%
\mathbf{E}\left( \nu \right) >1$. We are back to the discussion of
Subsection $2.1$ with $\left\{ X_{n}^{\overline{e}}\right\} $ either going
to $\infty $ or to a limiting rv depending on the tails of $\nu _{\overline{e%
}}$.

The only two cases that really matter are thus the case developed in
Subsection $2.1$ and case $\left( a\right) $ with state $0$ absorbing wp $1$%
, which was dealt with. We will therefore only consider the remaining first
case when $\left\{ X_{n}\right\} $ has state-space $\Bbb{N}.$

\section{Large $i$ estimates of $m_{i}\mathbf{=}\max_{j=1,...,i}\nu _{j}$}

We will use ideas stemming from limit laws for maxima of a large sample of
iid rvs in the continuum to give large $i$ estimates of $m_{i}\mathbf{=}%
\max_{j=1,...,i}\nu _{j},$ \cite{EKM}.

Let $X>0$ be some real-valued rv with density and no atom at $0$. Suppose $X$
has a finite mean $\mathbf{E}\left( X\right) $. Let $\overline{F}_{X}\left(
x\right) =\mathbf{P}\left( X>x\right) ,$ $x>0,$ be its complementary
probability distribution function (pdf). Define the law of some
integral-valued rv $\nu \in \Bbb{N}$ by: 
\begin{equation}
\mathbf{P}\left( \nu >j\right) =\mathbf{P}\left( X>j\right) ,\text{ }%
j=0,1,...  \label{1}
\end{equation}
Let $\overline{F}\left( j\right) =\mathbf{P}\left( \nu >j\right) $, $%
j=0,1,...$. With $\mathbf{E}\left( X\right) =\int_{0}^{\infty }\mathbf{P}%
\left( X>x\right) dx$ and $\mathbf{E}\left( \nu \right) =\sum_{j\geq 0}%
\overline{F}\left( j\right) $ we have $\mathbf{E}\left( \nu \right) -1<%
\mathbf{E}\left( X\right) <\mathbf{E}\left( \nu \right) $. This suggests
that if $\mathbf{E}\left( X\right) $ is large, $\mathbf{E}\left( \nu \right) 
$ is very close to $\mathbf{E}\left( X\right) .$\newline

\subsection{Maxima of a large sample of iid rvs in the continuum}

Let $M_{i}=\max \left( X_{1},...,X_{i}\right) $ with $\left( X_{i}\right)
_{i\geq 1}$ iid with $X_{1}\overset{d}{=}X.$

Two cases arise:

$\left( i\right) $ Von Mises case: With $a\left( x\right) >0,$ absolutely
continuous (with respect to Lebesgue measure) with density $a^{\prime
}\left( x\right) $ having $\lim a^{\prime }\left( x\right) =0$ as $%
x\rightarrow \infty $, consider 
\begin{equation*}
\mathbf{P}\left( X>x\right) =c\exp \left[ -\int^{x}\frac{dz}{a\left(
z\right) }\right] \text{, }c>0.
\end{equation*}
Then $a\left( x\right) =\mathbf{E}\left( X-x\mid X>x\right) $ is the mean
excess function with $a\left( x\right) /x\rightarrow 0$ as $x\rightarrow
\infty .$

Define $d_{i}$ by $\overline{F}_{X}\left( c_{i}\right) =1/i$ and $d_{i}$ by $%
c_{i}=a\left( c_{i}\right) .$ We have 
\begin{equation*}
d_{i}^{-1}\left( M_{i}-c_{i}\right) \overset{d}{\rightarrow }G\text{ as }%
i\rightarrow \infty ,
\end{equation*}
where $G$ has a Gumbel distribution $\mathbf{P}\left( G\leq x\right)
=e^{-e^{-x}}$, $x$ real. The sequence $c_{i}$ is increasing with $i$ with $%
c_{i}/i\rightarrow 0$ so at sublinear rate.

With $\gamma $ the Euler constant, it then holds that 
\begin{equation*}
d_{i}^{-1}\left( \mathbf{E}\left( M_{i}\right) -c_{i}\right) \rightarrow 
\mathbf{E}\left( G\right) =\gamma \text{ as }i\rightarrow \infty ,
\end{equation*}
so when $i$ gets large $\mathbf{E}\left( M_{i}\right) \sim c_{i}.$

$\left( ii\right) $ With $\alpha >0,$ suppose 
\begin{equation*}
\mathbf{P}\left( X>x\right) =x^{-\alpha }L\left( x\right) ,
\end{equation*}
where $L\left( x\right) $ is some slowly varying function at $\infty $, with 
$L\left( tx\right) /L\left( x\right) \rightarrow 1$ as $x\rightarrow \infty $
for all $t>0$. Defining $c_{i}$ by $\overline{F}_{X}\left( c_{i}\right) =1/i$
we have 
\begin{equation*}
c_{i}^{-1}M_{i}\overset{d}{\rightarrow }F\text{ as }i\rightarrow \infty ,
\end{equation*}
where $F$ has a Fr\'{e}chet distribution $\mathbf{P}\left( F\leq x\right)
=e^{-x^{-\alpha }}$, $x>0$ with $\mathbf{E}\left( F\right) =\Gamma \left(
1-1/\alpha \right) $ if $\alpha >1$, $=\infty $ if $\alpha \in \left(
0,1\right] $.

If $\alpha >1,$ with $c_{i}=i^{1/\alpha }L_{1}\left( i\right) $ for some
other slowly varying function $L_{1}$, it holds that 
\begin{equation*}
c_{i}^{-1}\mathbf{E}\left( M_{i}\right) \rightarrow \mathbf{E}\left(
F\right) =\Gamma \left( 1-1/\alpha \right) \text{ as }i\rightarrow \infty ,
\end{equation*}
so when $i$ gets large $\mathbf{E}\left( M_{i}\right) \sim \Gamma \left(
1-1/\alpha \right) c_{i}.$ And the sequence $c_{i}$ is increasing also at
sublinear rate.\newline

\subsection{Maxima of a large sample of discrete iid rvs: large $i$
estimation of $m_{i}$}

Let $m_{i}=\max \left( \nu _{1},...,\nu _{i}\right) $ with $\left( \nu
_{i}\right) _{i\geq 1}$ iid with $\nu _{1}\overset{d}{=}\nu $ and $\nu $'s
law given from $X$'s law as before$.$

Let $c_{i}$ be defined by $\overline{F}\left( c_{i}\right) =1/i$. In
general, it is not true that, upon scaling $m_{i}$, there is a proper weak
limit for this scaled rv, because in general, in the discrete setting $%
\overline{F}\left( j\right) /\overline{F}\left( j-1\right) \nrightarrow 1$ ($%
\mathbf{P}\left( \nu =j\right) /\mathbf{P}\left( \nu >j-1\right)
\nrightarrow 0$) as $j\rightarrow \infty $. All that can be said is that $%
m_{i}-c_{i}$ is tight (or bounded in probability), with $m_{i}/c_{i}%
\rightarrow 1$ in probability as $i\rightarrow \infty $. Also, if we are
interested in 
\begin{equation*}
\mathbf{E}\left( m_{i}\right) =\sum_{j\geq 0}\left( 1-F\left( j\right)
^{i}\right) ,
\end{equation*}
using the latter argument, for large $i$, $\mathbf{E}\left( m_{i}\right) $
and $\mathbf{E}\left( M_{i}\right) $ are of the same order of magnitude.

Ergodic case from (\ref{L1})$:$ With $c_{i}$ defined by $\overline{F}%
_{X}\left( c_{i}\right) =1/i$ in $\left( i\right) $ the Von Mises case or $%
\left( ii\right) $ when $\mathbf{P}\left( X>x\right) =x^{-\alpha }L\left(
x\right) $ and $\alpha >1$ in the domain of attraction of the Fr\'{e}chet$%
\left( \alpha \right) $ law

\begin{eqnarray*}
\left( i\right) \text{ }\mathbf{E}\left( m_{i}\right) &\sim &c_{i} \\
\left( ii\right) \text{ }\mathbf{E}\left( m_{i}\right) &\sim &\Gamma \left(
1-1/\alpha \right) c_{i}=\Gamma \left( 1-1/\alpha \right) i^{1/\alpha
}L_{1}\left( i\right) .
\end{eqnarray*}
In all these cases, $\mathbf{E}\left( m_{i}\right) $ grows at sublinear rate
as $i$ gets large.\newline

Under the above assumptions on the law of $\nu $, there is thus an integer $%
I $ such that 
\begin{eqnarray*}
\mathbf{E}\left( X_{n+1}\mid X_{n}=i\right) &\leq &i-1\text{ for all }i\geq I
\\
\mathbf{E}\left( X_{n+1}\mid X_{n}=i\right) &<&\infty \text{ for all }i\text{
for which }\mathbf{E}\left( X_{n+1}\mid X_{n}=i\right) >i-1,
\end{eqnarray*}
which by Foster theorem implies that $\left\{ X_{n}\right\} $ is ergodic, 
\cite{Fost}. The limit law of the MC is the unique integrable solution to
the corresponding functional equation for $X_{\infty }.$\newline

Transient case from (\ref{L2}): If $\nu $ is in the domain of attraction of
the Fr\'{e}chet law with $\alpha \in \left( 0,1\right) $, $\mathbf{E}\left(
m_{i}\right) =\mathbf{E}\left( X_{n+1}\mid X_{n}=i\right) $ grows at a
superlinear rate which leads to a transience case ($X_{n}\rightarrow \infty $
wp $1$, as $n\rightarrow \infty $). Such $\nu $s have infinite mean. If $%
\alpha =1$, the process is transient (positive recurrent) if $\lim \sup_{i}i%
\mathbf{P}\left( \nu >i\right) >e^{-\gamma }$ (respectively $<e^{-\gamma }$%
). Whenever the tails of $\nu $ satisfy any one of these conditions, $%
\mathbf{E}\left( \nu \right) =\infty .$ This shows that transience of $%
\left\{ X_{n}\right\} $ does not necessarily mean $\mathbf{E}\left( \nu
\right) =\infty $.\newline

\emph{Example:} If $\alpha =1$, there are positive recurrent examples for
which $\mathbf{E}\left( \nu \right) =\infty $, for instance those obtained
from 
\begin{equation*}
\mathbf{P}\left( \nu >i\right) =\frac{1}{i\log ^{\beta }\left( 1+i\right) }%
\text{ with }0<\beta <1\text{, }
\end{equation*}
with $L\left( x\right) =\log ^{\beta }\left( 1+x\right) $ slowly varying at $%
\infty $. $\diamondsuit $

\section{General approach to find solutions of the functional equation}

In the ergodic case from (\ref{L1}), the invariant probability measure $\pi
\left( j\right) :=\mathbf{P}\left( X_{\infty }=j\right) $ solves 
\begin{equation*}
\mathbf{\pi }^{\prime }=\mathbf{\pi }^{\prime }P.
\end{equation*}
However, here, the pdfs of $\nu >0$ and $X_{\infty }>0$ are related by the
functional equation: 
\begin{equation}
F_{\infty }\left( j\right) =\Phi _{\infty }\left( F\left( j\right) \right) ,
\label{FE}
\end{equation}
and we shall give many examples of explicit pairs $\left( F_{\infty }\left(
j\right) ,F\left( j\right) \right) $ solving it. As indicated above, it
participates to the general program of finding the branching number $\nu $
of the Lamperti mbp realizing an initial target guess of the limiting rv $%
X_{\infty }$. The rv $X_{\infty }$ is taken from the classical (shifted) set
of probability mass functions (pmfs) supported by the integers. We will then
compute explicitly the law of $\nu $ corresponding to classical target pmfs
such as geometric, Sibuya, Poisson. The obtained distributions are far from
classical and somewhat surprising.\newline

\emph{Remark:} With the idea of spanning trees in the background, there
exists a determinantal Kirchoff formula stating that, \cite{Pit}: 
\begin{equation*}
\pi \left( j\right) =\det \left[ \left( I-P\right) ^{\left( j,j\right)
}\right] ,
\end{equation*}
where $\left( I-P\right) ^{\left( j,j\right) }$ is the Laplacian matrix $I-P$
to which row $j$ and column $j$ have been removed. In view of the expression
(\ref{P1} with $i,j\geq 1$) of the Lamperti matrix $P$, the Kirchoff formula
shows that the computation of $\mathbf{\pi }$ from $P$ (and so from $F$) is
not, in principle, a simple matter. Our approach being to find $F$ (and so $%
P $), starting from the knowledge of $\mathbf{\pi },$ this leads in return
to non trivial determinantal identities. $\diamondsuit $\newline

- \textbf{Lagrange inversion formula:}

In the sequel, we shall denote by:

$\left( n\right) _{k}=n\left( n-1\right) ...\left( n-k+1\right) $ and $%
\left[ n\right] _{k}=n\left( n+1\right) ...\left( n+k-1\right) $ the falling
and rising factorials (of order $k$) of $n$.

Take $\Phi _{\infty }\left( z\right) =z\Psi _{\infty }\left( z\right) $ for
some new (given) pgf $\Psi _{\infty }$ obeying $\Psi _{\infty }\left(
0\right) \neq 0.$ The pgf $\Psi _{\infty }$ is the one of $X_{\infty }-1$.
Apply Lagrange inversion formula to solve $z\Psi _{\infty }\left( z\right)
=u $. It gives the inverse $\Phi _{\infty }^{-1}\left( z\right) $ of $\Phi
_{\infty }\left( z\right) $ as a power series in $z,$ with 
\begin{equation*}
\varphi _{n}:=\left[ z^{n}\right] \Phi _{\infty }^{-1}\left( z\right) =\frac{%
1}{n}\left[ z^{n-1}\right] \Psi _{\infty }\left( z\right) ^{-n}.
\end{equation*}
Then 
\begin{equation*}
F\left( j\right) =\Phi _{\infty }^{-1}\left( F_{\infty }\left( j\right)
\right) =\sum_{n\geq 1}\varphi _{n}F_{\infty }\left( j\right) ^{n}
\end{equation*}
gives the $F\left( j\right) $ consistent with the original choice of $%
F_{\infty }\left( j\right) $. With $B_{n,k}\left( x_{1},x_{2},...\right) $
(respectively $\widehat{B}_{n,k}\left( x_{1},x_{2},...\right) $) the
exponential (respectively ordinary) Bell polynomials in the indeterminates $%
\left( x_{1},x_{2},...\right) $, obeying $B_{n,k}\left(
x_{1},x_{2},...\right) =0$ if $k>n$ and $B_{n,0}\left(
x_{1},x_{2},...\right) =\delta _{n,0},$ we have in principle (\cite{Com}, p. 
$161$) ($x_{k}=k!\pi \left( k+1\right) $) 
\begin{eqnarray*}
\varphi _{n} &=&\frac{1}{n!}\sum_{k=0}^{n-1}\left( -n\right) _{k}\pi \left(
1\right) ^{-\left( n+k\right) }B_{n-1,k}\left( \frac{2!\pi \left( 2\right) }{%
2},\frac{3!\pi \left( 3\right) }{3},...\right) \\
&=&\frac{1}{n}\sum_{k=0}^{n-1}\frac{\left( -n\right) _{k}}{k!}\pi \left(
1\right) ^{-\left( n+k\right) }\widehat{B}_{n-1,k}\left( \pi \left( 2\right)
,\pi \left( 3\right) ,...\right) \\
&=&\frac{\pi \left( 1\right) ^{-n}}{n}\sum_{k=0}^{n-1}\frac{\left( -n\right)
_{k}}{k!}\widehat{B}_{n-1,k}\left( \frac{\pi \left( 2\right) }{\pi \left(
1\right) },\frac{\pi \left( 3\right) }{\pi \left( 1\right) },...\right)
\end{eqnarray*}
\begin{eqnarray*}
&=&\frac{1}{n!}\sum_{k=0}^{n-1}\left( -1\right) ^{k}\pi \left( 1\right)
^{-\left( n+k\right) }B_{n+k-1,k}\left( 0,2!\pi \left( 2\right) ,3!\pi
\left( 3\right) ,...\right) \\
&=&\frac{\pi \left( 1\right) ^{-n}}{n}\sum_{k=0}^{n-1}\left( -1\right) ^{k}%
\binom{n+k-1}{k}\widehat{B}_{n+k-1,k}\left( 0,\frac{\pi \left( 2\right) }{%
\pi \left( 1\right) },\frac{\pi \left( 3\right) }{\pi \left( 1\right) }%
,...\right) .
\end{eqnarray*}
Owing to $\left( -n\right) _{k}=\left( -1\right) ^{k}\left[ n\right] _{k}$
and (see \cite{Com}, p. $145$) 
\begin{eqnarray*}
B_{n-1,k}\left( x_{1},x_{2},...\right) &=&\left( n-1\right)
!\sum^{*}\prod_{m\geq 1}\frac{1}{k_{m}!}\left( \frac{x_{m}}{m!}\right)
^{k_{m}}, \\
\widehat{B}_{n-1,k}\left( x_{1},x_{2},...\right) &=&k!\sum^{*}\prod_{m\geq 1}%
\frac{x_{m}^{k_{m}}}{k_{m}!}
\end{eqnarray*}
where the star sum runs over $k_{m}\geq 0$, obeying $\sum_{m\geq 1}k_{m}=k$
and $\sum_{m\geq 1}mk_{m}=n-1$, we have equivalently $\varphi _{1}=1/\pi
\left( 1\right) $ and if $n\geq 2$%
\begin{equation*}
\varphi _{n}=\frac{\pi \left( 1\right) ^{-n}}{n}\sum_{k=1}^{n-1}\left(
-1\right) ^{k}\left[ n\right] _{k}C_{n-1,k},
\end{equation*}
where, with $C_{n-1,0}=\delta _{n-1,0},$%
\begin{equation}
C_{n-1,k}=\sum^{*}\prod_{m\geq 1}\frac{\left( \pi \left( m+1\right) /\pi
\left( 1\right) \right) ^{k_{m}}}{k_{m}!}.  \label{C}
\end{equation}
To summarize, we obtained the pdf $F$ of $\nu $ corresponding to $\pi \left(
j\right) :=\mathbf{P}\left( X_{\infty }=j\right) ,$ solving (\ref{FE}), as

\begin{proposition}
The mapping $X_{\infty }\rightarrow \nu $ is one-to-one. With $C_{n-1,k}$
given by (\ref{C}) and $h_{n}=\frac{1}{n}\sum_{k=1}^{n-1}\left( -1\right)
^{k}\left[ n\right] _{k}C_{n-1,k}$ ($h_{1}=1$), 
\begin{equation}
F\left( j\right) =\Phi _{\infty }^{-1}\left( F_{\infty }\left( j\right)
\right) =\sum_{n\geq 1}h_{n}\cdot \left( F_{\infty }\left( j\right) /\pi
\left( 1\right) \right) ^{n}  \label{FES}
\end{equation}
is the cumulated mass function of $\nu $ corresponding to any given $\mathbf{%
\pi }$.
\end{proposition}

The obtained expression (\ref{FES}) only depends on the ratio $F_{\infty
}\left( j\right) /F_{\infty }\left( 1\right) $. Note that $\Phi _{\infty
}^{-1}\left( z\right) $ is increasing from $z=0$ to $z=1$ and concave. From
the fact that it is increasing, we conclude that if $F_{\infty }\left(
j\right) $ is a pdf, then so is $F\left( j\right) $. From the concavity, we
conclude $F\left( j\right) \geq F_{\infty }\left( j\right) $ for all $j$ (as
already mentioned, $X_{\infty }$ is stochastically larger than $\nu $).
While proceeding in this way, we observe that, given we first fix the law of 
$X_{\infty }$, the one of the corresponding $\nu $ follows.\newline

Suppose we were able to find a suitable pair of pdfs $\left( F\left(
j\right) ,F_{\infty }\left( j\right) \right) $ by the Lagrange inversion
formula. Then, with $F_{0}\left( j\right) =\mathbf{1}\left( j\leq 1\right) $
($X_{0}\overset{d}{\sim }\delta _{1}$) and $\Phi _{0}\left( z\right) =z$, $%
F_{1}\left( j\right) =\Phi _{0}\left( F\left( j\right) \right) =F\left(
j\right) $ is a pdf, the one of $\nu $. Let $\Phi _{1}\left( z\right)
=\sum_{j\geq 1}z^{j}\left( F_{1}\left( j\right) -F_{1}\left( j-1\right)
\right) $ be the pgf of $X_{1}\overset{d}{=}\nu $. Next, $F_{2}\left(
j\right) =\Phi _{1}\left( F\left( j\right) \right) $ is a pdf because $\Phi
_{1}$ is monotone increasing obeying $\Phi _{1}\left( 0\right) =0$, $\Phi
_{1}\left( 1\right) =1$. By recurrence $F_{n+1}\left( j\right) =\Phi
_{n}\left( F\left( j\right) \right) $ is the pdf of some rv $X_{n+1}$
obtained from the one of $X_{n}$ and $F_{n}\left( j\right) \rightarrow
F_{\infty }\left( j\right) $ solution to $F_{\infty }\left( j\right) =\Phi
_{\infty }\left( F\left( j\right) \right) $.\newline

We shall deal with special cases of $X_{\infty }.$

- \textbf{Infinite divisibility:} suppose $\Psi _{\infty }\left( z\right) $
is the pgf of an infinitely divisible (ID) rv (meaning $X_{\infty }-1$ is
ID). Then, as a compound Poisson rv, 
\begin{equation*}
\Psi _{\infty }\left( z\right) =e^{-\lambda \left( 1-h\left( z\right)
\right) },
\end{equation*}
for some rate $\lambda >0$ and pgf $h\left( z\right) $ obeying $h\left(
0\right) =0.$ If $\mathbf{P}\left( X_{\infty }=j+1\right) =\left[
z^{j}\right] \Psi _{\infty }\left( z\right) =\pi _{\lambda }\left(
j+1\right) $ is a known simple function of $\lambda $, then $\left[
z^{j}\right] \Psi _{\infty }\left( z\right) ^{-n}$ is readily obtained as $%
\pi _{-n\lambda }\left( j+1\right) $, a useful identity to get the $h_{n}$
in (\ref{FES}) and so $F$ from $F_{\infty }$. \newline

- \textbf{Complete monotonicity:} Suppose $\overline{F}_{\infty }\left(
j\right) $ defines a in $\left[ 0,1\right] $-valued completely monotone
sequence of complementary pdfs, meaning 
\begin{eqnarray*}
\left( -1\right) ^{k}\Delta ^{\left( k\right) }\overline{F}_{\infty }\left(
j\right) &\geq &0\text{ for all }j,k\geq 0,\text{ equivalently} \\
\left( -1\right) ^{k}\Delta ^{\left( k\right) }\mathbf{P}\left( X_{\infty
}=j\right) &\geq &0\text{ for all }j\geq 1,k\geq 0,
\end{eqnarray*}
where $\Delta :$ $\Delta h\left( j\right) =h\left( j+1\right) -h\left(
j\right) $ is the right-shift operator and $\Delta ^{\left( k\right) }$ its $%
k-$th iterate. Note $\mathbf{P}\left( X_{\infty }=j\right) =\Delta F_{\infty
}\left( j-1\right) =-\Delta \overline{F}_{\infty }\left( j-1\right) $. By
Hausdorff representation theorem, $\overline{F}_{\infty }\left( j\right) $
is completely monotone (CM) if and only if 
\begin{equation*}
\overline{F}_{\infty }\left( j\right) =\int_{0}^{1}u^{j}\lambda \left(
du\right) ,
\end{equation*}
for some probability measure $\lambda \left( du\right) $ on $\left[
0,1\right] .$

Equivalently, with $\Phi _{\infty }\left( z\right) =\sum_{j\geq 1}z^{j}%
\mathbf{P}\left( X_{\infty }=j\right) $ the pgf of $X_{\infty }$, 
\begin{equation*}
\sum_{j\geq 0}z^{j}\overline{F}_{\infty }\left( j\right) =\frac{1-\Phi
_{\infty }\left( z\right) }{1-z}=\int_{0}^{1}\frac{1}{1-zu}\lambda \left(
du\right) ,
\end{equation*}
as a Stieltjes transform of $\lambda \left( du\right) .$ Note that, with $U%
\overset{d}{\sim }\lambda \left( du\right) $%
\begin{equation*}
\Phi _{\infty }\left( z\right) =z\int_{0}^{1}\frac{1-u}{1-zu}\lambda \left(
du\right) =\mathbf{E}\left( \frac{z\left( 1-U\right) }{1-zU}\right) ,
\end{equation*}
showing that $X_{\infty }-1$, with pgf $\Psi _{\infty }\left( z\right)
=z^{-1}\Phi _{\infty }\left( z\right) =\mathbf{E}\left( \frac{1-U}{1-zU}%
\right) $, is a $\lambda -$mixture of a shifted geometric rv, so that $%
X_{\infty }-1$ is log-convex and infinitely divisible. As noted in \cite
{Gupta}, log-convex (log-concave) pmfs are decreasing (increasing) failure
rate monotone, say DFR (IFR), meaning $\Delta r_{j}$ decreasing (increasing)
where $r_{j}=\pi \left( j\right) /\overline{F}_{\infty }\left( j-1\right) =%
\mathbf{P}\left( X_{\infty }=j\right) /\mathbf{P}\left( X_{\infty }\geq
j\right) $ is a discrete failure `rate'.

\subsection{Explicit examples of $\left( \nu ,X_{\infty }\right) $ with
support $\left\{ 1,...,\infty \right\} $.}

In some cases, the computation of the pair $\left( F\left( j\right)
,F_{\infty }\left( j\right) \right) $ is obtained as a simple expression.

- \textbf{Geometric} example: $X_{\infty }\overset{d}{\sim }$geom$\left(
p\right) $

\begin{proposition}
Suppose $X_{\infty }\overset{d}{\sim }$geom$\left( p\right) $, so with $%
F_{\infty }\left( j\right) =1-q^{j}$. The sequence $\overline{F}_{\infty
}\left( j\right) $ is of course CM as a result of 
\begin{equation*}
\overline{F}_{\infty }\left( j\right) =q^{j}=\int_{0}^{1}u^{j}\lambda \left(
du\right) ,\text{ with }\lambda \left( du\right) =\delta _{q}\left(
du\right) ,
\end{equation*}
so $X_{\infty }-1$ is ID.

$\left( i\right) $ The solution to (\ref{FE}) is: 
\begin{equation*}
\mathbf{P}\left( \nu \leq j\right) =F\left( j\right) =\frac{1-q^{j}}{%
1-q^{j+1}}\text{, }j=1,2,...
\end{equation*}

$\left( ii\right) $ The sequence $\overline{F}\left( j\right) $ is CM and so 
$\nu -1$ is ID. The distribution $F\left( j\right) $ has decreasing failure
rate (DFR).

$\left( iii\right) $ There are two ways to generate the corresponding
branching number $\nu :$%
\begin{equation*}
\left( iii-a\right) :\text{ }\nu =\inf \left( i\geq 1:\mathcal{B}_{i}\left(
\alpha _{i}\right) =1\right) ,
\end{equation*}
where $\left( \mathcal{B}_{i}\left( \alpha _{i}\right) ;\text{ }i\geq
1\right) $ is an independent sequence of Bernoulli rvs with success
parameter $\alpha _{i}=1/\left( 1+q+...+q^{i}\right) .$ Or: 
\begin{equation*}
\left( iii-b\right) :\text{ }\nu =\max_{i=1,...,G}\xi _{i}
\end{equation*}
where $G\overset{d}{\sim }$geom$\left( p\right) $ independent of the iid
sequence $\left( \xi _{i}\text{, }i\geq 1\right) $ with $\xi _{1}\overset{d}{%
\sim }$geom$\left( p\right) .$

$\left( iv\right) $ The tails of both $\left( \nu ,X_{\infty }\right) $ are
geometric with: $\mathbf{P}\left( \nu >j\right) /\mathbf{P}\left( X_{\infty
}>j\right) \rightarrow p<1$.
\end{proposition}

\emph{Proof:}

$\left( i\right) $ We have 
\begin{equation*}
\Phi _{\infty }\left( z\right) =\mathbf{E}\left( z^{X_{\infty }}\right) =%
\frac{pz}{1-qz}\text{ and so}
\end{equation*}
\begin{equation*}
\Phi _{\infty }\left( F\left( j\right) \right) =\frac{p\frac{1-q^{j}}{%
1-q^{j+1}}}{1-q\left( \frac{1-q^{j}}{1-q^{j+1}}\right) }=1-q^{j}=\mathbf{P}%
\left( X_{\infty }\leq j\right) =F_{\infty }\left( j\right) .
\end{equation*}

$\left( ii\right) $ With $\lambda \left( du\right) =p\sum_{j\geq
1}q^{j-1}\delta _{q^{j}}$, a probability measure, 
\begin{equation*}
\overline{F}\left( j\right) =\frac{pq^{j}}{1-q^{j+1}}=\int_{0}^{1}u^{j}%
\lambda \left( du\right) .
\end{equation*}
The rv $\nu \geq 1$ has finite mean $\mathbf{E}\left( \nu \right)
=\sum_{j\geq 0}\mathbf{P}\left( \nu >j\right) =1+p\sum_{j\geq 1}q^{j}/\left(
1-q^{j+1}\right) <\infty $ and $\left\{ X_{n}\right\} $ is recurrent
positive.

For $j\geq 1$, we have 
\begin{equation*}
\frac{\mathbf{P}\left( \nu =j\right) }{\mathbf{P}\left( \nu >j\right) }=%
\frac{p}{q}\frac{1}{1-q^{j}}
\end{equation*}
which is decreasing with $j$.

$\left( iii\right) $ The first statement $\left( iii-a\right) $ results
from: $\mathbf{P}\left( \nu >i\right) =\prod_{j=1}^{i}\left( 1-\alpha
_{j}\right) =\frac{pq^{i}}{1-q^{i+1}}.$

$\left( iii-b\right) $ results from 
\begin{equation*}
\mathbf{P}\left( \max_{i=1,...,G}\xi _{i}>i\right) =\sum_{k\geq
1}pq^{k-1}q^{ik}=\frac{pq^{i}}{1-q^{i+1}}.
\end{equation*}

$\left( iv\right) $ The tails of $\nu $ are given by $\mathbf{P}\left( \nu
>j\right) =1-\Phi _{\infty }^{-1}\left( \mathbf{P}\left( X_{\infty }\leq
j\right) \right) =1-\frac{z}{p+qz}\mid _{1-p^{j}}\sim pq^{j}$ (for large $j$%
) with $\mathbf{P}\left( \nu >j\right) /\mathbf{P}\left( X_{\infty
}>j\right) \rightarrow p<1$. For this model, $\nu $ and $X_{\infty }$ are
geometric (power-law) and tail-equivalent but the tails of $\nu $ are
thinner than the ones of $X_{\infty }.$ $\Box $\newline

Related examples to the geometric one ($X_{\infty }$ having dominant
geometric tails with an algebraic prefactor):

- Suppose $X_{\infty }\overset{d}{\sim }$\textbf{negative-binomial ,}%
conditioned to be positive\textbf{: }With\textbf{\ }$\left[ \alpha \right]
_{k}=\Gamma \left( \alpha +k\right) /\Gamma \left( \alpha \right) $, $\alpha
>0$, suppose 
\begin{equation*}
\Phi _{\infty }\left( z\right) =\mathbf{E}\left( z^{X_{\infty }}\right) =%
\frac{\left( \frac{p}{1-qz}\right) ^{\alpha }-p^{\alpha }}{1-p^{\alpha }}
\end{equation*}
is the pgf of a negative-binomial\textbf{\ }rv\textbf{, }conditioned to be
positive. Then, by direct inversion 
\begin{equation*}
F\left( j\right) =\Phi _{\infty }^{-1}\left( F_{\infty }\left( j\right)
\right) =\frac{1-\left( 1+\sum_{k=1}^{j}\frac{\left[ \alpha \right] _{k}}{k!}%
q^{k}\right) ^{-1/\alpha }}{q},
\end{equation*}
which defines the pdf of $\nu $. In this case, $\Psi _{\infty }\left(
z\right) =z^{-1}\Phi _{\infty }\left( z\right) $ is not the pgf of an ID rv.
Plugging $\alpha =1$ gives back the latter geometric case. The tails of $%
X_{\infty }$ goes, up to a constant prefactor, like $j^{\alpha -1}q^{j}$.

- Suppose $X_{\infty }\overset{d}{\sim }$shifted \textbf{negative-bin (}$%
\Psi _{\infty }\left( z\right) $ now is the pgf of the ID rv $X_{\infty }-1$%
): then 
\begin{eqnarray*}
\Phi _{\infty }\left( z\right) &=&\mathbf{E}\left( z^{X_{\infty }}\right)
=z\left( \frac{p}{1-qz}\right) ^{\alpha }\text{ and }\Psi _{\infty }\left(
z\right) =p^{\alpha }\left( 1-qz\right) ^{-\alpha } \\
\left[ z^{j}\right] \Psi _{\infty }\left( z\right) &=&p^{\alpha }\frac{%
\left[ \alpha \right] _{j}}{j!}q^{j}\Rightarrow \left[ z^{j}\right] \Psi
_{\infty }^{-n}\left( z\right) =\frac{\left[ -n\alpha \right] _{j}}{j!}%
p^{-n\alpha }q^{j} \\
\frac{1}{n}\left[ z^{n-1}\right] \Psi _{\infty }\left( z\right) ^{-n} &=&%
\frac{\left[ -n\alpha \right] _{n-1}}{n!}p^{-n\alpha }q^{n-1} \\
F_{\infty }\left( j\right) &=&\sum_{k=1}^{j}\mathbf{P}\left( X_{\infty
}=k\right) =p^{\alpha }\sum_{k=1}^{j}\frac{\left[ \alpha \right] _{k-1}}{%
\left( k-1\right) !}q^{k-1} \\
F\left( j\right) &=&\Phi _{\infty }^{-1}\left( F_{\infty }\left( j\right)
\right) =\sum_{n\geq 1}F_{\infty }\left( j\right) ^{n}\frac{1}{n}\left[
z^{n-1}\right] \Psi _{\infty }\left( z\right) ^{-n}.
\end{eqnarray*}
The negative-binomial distribution with pgf $\Psi _{\infty }\left( z\right)
=p^{\alpha }\left( 1-qz\right) ^{-\alpha }$ is CM (and so log-convex, ID and
DFR) only when $\alpha \leq 1$. When $\alpha \geq 1$ it is log-concave, ID
and IFR.

- $X_{\infty }\overset{d}{\sim }$\textbf{Fisher log-series}. With $p\in
\left( 0,1\right) $ and $c=-\log \left( 1-p\right) $, suppose 
\begin{eqnarray*}
\Phi _{\infty }\left( z\right) &=&\mathbf{E}\left( z^{X_{\infty }}\right)
=-c^{-1}\log \left( 1-pz\right) =:z\Psi _{\infty }\left( z\right) \\
\mathbf{P}\left( X_{\infty }=k\right) &=&c^{-1}p^{k}/k\text{, }k\geq 1\text{
and }\mathbf{P}\left( X_{\infty }\leq j\right) =c^{-1}\sum_{k=1}^{j}p^{k}/k%
\text{ }
\end{eqnarray*}
involving a truncated logarithm. We have 
\begin{eqnarray*}
\Phi _{\infty }^{-1}\left( z\right) &=&p^{-1}\left( 1-e^{-cz}\right) \text{
and } \\
F\left( j\right) &=&\Phi _{\infty }^{-1}\left( F_{\infty }\left( j\right)
\right) =p^{-1}\left( 1-e^{-\sum_{k=1}^{j}p^{k}/k}\right) .
\end{eqnarray*}
For all $j\geq 1$, we have by construction 
\begin{equation*}
F\left( j\right) >F_{\infty }\left( j\right) =c^{-1}\sum_{k=1}^{j}p^{k}/k.
\end{equation*}
The tails of $X_{\infty }$ goes, up to a constant prefactor, like $%
j^{-1}p^{j}$. In addition, 
\begin{eqnarray*}
\mathbf{P}\left( X_{\infty }=i\right) &=&\int_{0}^{1}u^{i-1}\mu \left(
du\right) \text{ where }\mu \left( du\right) =c^{-1}1_{u\in \left(
0,p\right) }du \\
\mathbf{P}\left( X_{\infty }>i\right) &=&\int_{0}^{1}u^{i}\lambda \left(
du\right) \text{ where }\lambda \left( du\right) =c^{-1}\left( 1-u\right)
^{-1}1_{u\in \left( 0,p\right) }du,
\end{eqnarray*}
and both $X_{\infty }$ and $\nu $ are CM. \newline

Let us now look at situations when $X_{\infty }$ has heavy (algebraic) tails
with index $\alpha >0$:

- The power-law \textbf{Sibuya} example, \cite{Sibu}.

\begin{proposition}
With $\alpha \in \left( 0,1\right) $, suppose $X_{\infty }\overset{d}{\sim }$%
Sibuya$\left( \alpha \right) $, that is:

$\Phi _{\infty }\left( z\right) =\mathbf{E}\left( z^{X_{\infty }}\right)
=1-\left( 1-z\right) ^{\alpha },$ with $\mathbf{P}\left( X_{\infty
}=j\right) =\pi \left( j\right) =\alpha \left[ 1-\alpha \right] _{j-1}/j!$.
Then:

$\left( i\right) $ The sequence $\pi \left( j\right) $ is CM, so log-convex,
DFR and $X_{\infty }-1$ is ID.

$\left( ii\right) $ 
\begin{equation}
X_{\infty }=\inf \left( i\geq 1:\mathcal{B}_{i}\left( \alpha _{i}\right)
=1\right) ,  \label{ber}
\end{equation}
where $\left( \mathcal{B}_{i}\left( \alpha _{i}\right) \right) _{i\geq 1}$
is a sequence of independent Bernoulli rvs obeying $\mathbf{P}\left( 
\mathcal{B}_{i}\left( \alpha _{i}\right) =1\right) =\alpha /i.$

$\left( iii\right) $ The solution to (\ref{FE}) is: 
\begin{equation*}
\mathbf{P}\left( \nu \leq i\right) =1-\left( 1-\alpha \sum_{j=1}^{i}\left[
1-\alpha \right] _{j-1}/j!\right) ^{1/\alpha },
\end{equation*}

$\left( iv\right) $ Both $\left( X_{\infty },\nu \right) $ have algebraic
(power-law) tails, but with tail index $\alpha $ and $1$ respectively.

$\left( v\right) $ We have 
\begin{equation*}
\mathbf{P}\left( \nu >i\right) \sim \frac{1}{\Gamma \left( 1-\alpha \right)
^{1/\alpha }}i^{-1}\text{ as }j\rightarrow \infty 
\end{equation*}
and $1/\Gamma \left( 1-\alpha \right) ^{1/\alpha }<e^{-\gamma }.$ For all $%
\alpha \in \left( 0,1\right) $, the Lamperti chain generated by $\nu $ is
positive recurrent, with invariant probability measure $\mathbf{\pi }$.
\end{proposition}

\emph{Proof:} $\left( i\right) $ It can be checked that, with $\mu \left(
du\right) \overset{d}{\sim }$Beta$\left( 1-\alpha ,\alpha \right) $%
\begin{equation*}
\mathbf{P}\left( X_{\infty }=j\right) =\int_{0}^{1}u^{j}\mu \left( du\right)
.
\end{equation*}

$\left( ii\right) $ is obvious and a well-known property of Sibuya$\left(
\alpha \right) $ distributed rvs, \cite{Sibu}.

$\left( iii\right) $ We have $\Phi _{\infty }^{-1}\left( z\right) =1-\left(
1-z\right) ^{1/\alpha }$ and so

\begin{equation*}
\mathbf{P}\left( \nu \leq i\right) =1-\left( 1-\alpha \sum_{j=1}^{i}\left[
1-\alpha \right] _{j-1}/j!\right) ^{1/\alpha }.
\end{equation*}

$\left( iv\right) $ We have $\mathbf{P}\left( X_{\infty }=j\right) \sim 
\frac{\alpha }{\Gamma \left( 1-\alpha \right) }j^{-\left( \alpha +1\right) }$
and $\mathbf{P}\left( X_{\infty }>j\right) \sim \frac{1}{\Gamma \left(
1-\alpha \right) }j^{-\alpha }.$ Therefore $\mathbf{P}\left( \nu \leq
j\right) \sim \Phi _{\infty }^{-1}\left( \mathbf{P}\left( X_{\infty }\leq
j\right) \right) \sim 1-\mathbf{P}\left( X_{\infty }>j\right) ^{1/\alpha
}\sim 1-\frac{1}{\Gamma \left( 1-\alpha \right) ^{1/\alpha }}j^{-1}.$ And $%
\nu $ has lighter tails (of index $1$) than $X_{\infty }$ (of index $\alpha $%
)$.$ This is a concrete manifestation in the tails of the fact that $%
X_{\infty }$ is stochastically larger than $\nu $.

$\left( v\right) $ To decide whether or not $\nu $ belongs to the ergodic
family, (\ref{L1}), we need to compare $1/\Gamma \left( 1-\alpha \right)
^{1/\alpha }$ with $e^{-\gamma }$, $\gamma =-\Gamma ^{\prime }\left(
1\right) $ being the Euler constant. Indeed, based on Lamperti's criterion,
the chain is recurrent if $1/\Gamma \left( 1-\alpha \right) ^{1/\alpha
}<e^{-\gamma }$ or $\log \Gamma \left( 1-\alpha \right) /\alpha >\gamma $
for all $\alpha \in \left( 0,1\right) $. But this is always true because $%
\log \Gamma \left( 1-\alpha \right) /\alpha $ is an increasing function of $%
\alpha $ with $\log \Gamma \left( 1-\alpha \right) /\alpha \rightarrow
\gamma $ as $\alpha \rightarrow 0$ ($\log \Gamma \left( 1-\alpha \right)
\sim \log \left( 1-\alpha \Gamma ^{\prime }\left( 1\right) \right) \sim
-\alpha \Gamma ^{\prime }\left( 1\right) $)$.$ The critical upper bound $%
e^{-\gamma }$ for the coefficient $1/\Gamma \left( 1-\alpha \right)
^{1/\alpha }$ is attained for $\alpha \rightarrow 0.$ $\Box $\newline

Related examples to the Sibuya one with power-law tails are:

- \textbf{Pareto} ($\alpha >0$): Suppose $\mathbf{P}\left( X_{\infty
}>i\right) =\left( i+1\right) ^{-\alpha }$. Clearly, 
\begin{equation}
X_{\infty }=\inf \left( i\geq 1:\mathcal{B}_{i}\left( \alpha _{i}\right)
=1\right) ,
\end{equation}
where $\left( \mathcal{B}_{i}\left( \alpha _{i}\right) \right) _{i\geq 1}$
is a sequence of independent Bernoulli rvs obeying $\mathbf{P}\left( 
\mathcal{B}_{i}\left( \alpha _{i}\right) =1\right) =1-\left( 1+1/i\right)
^{-\alpha }$ where $\alpha >0$. Indeed, 
\begin{equation*}
\mathbf{P}\left( X_{\infty }>i\right) =\prod_{j=1}^{i}\left( 1-\alpha
_{j}\right) =\prod_{j=1}^{i}\left( 1+1/j\right) ^{-\alpha }=\left(
i+1\right) ^{-\alpha }.
\end{equation*}
We have $\mathbf{P}\left( X_{\infty }=i\right) =i^{-\alpha }-\left(
i+1\right) ^{-\alpha }=i^{-\alpha }\left( 1-\left( \left( i+1\right)
/i\right) ^{-\alpha }\right) \sim \alpha i^{-\left( \alpha +1\right) }$ and
so $\Phi _{\infty }\left( z\right) =\sum_{i\geq 1}z^{i}i^{-\alpha
}-\sum_{i\geq 1}z^{i}\left( i+1\right) ^{-\alpha }=1-z^{-1}\left( 1-z\right)
L_{\alpha }\left( z\right) =z\Psi _{\infty }\left( z\right) .$

When $\alpha \leq 1$, the polylog function $L_{\alpha }\left( z\right)
=\sum_{i\geq 1}z^{i}i^{-\alpha }$ is not defined at $z=1$ but $z\Psi
_{\infty }\left( z\right) =1-z^{-1}\left( 1-z\right) L_{\alpha }\left(
z\right) $ is a true pgf taking the value $1$ at $z=1$. Lagrange inversion
formula gives the power-series expansion of $\Phi _{\infty }^{-1}\left(
z\right) $ giving $\mathbf{P}\left( \nu \leq j\right) =\Phi _{\infty
}^{-1}\left( 1-\left( j+1\right) ^{-\alpha }\right) .$

The rv $X_{\infty }-1$ (with pgf $\Psi _{\infty }\left( z\right) $) is
infinitely divisible. Indeed, the polylogarithm can be expressed in terms of
the integral of the Bose-Einstein distribution

\begin{equation*}
L_{\alpha }\left( z\right) =\frac{1}{\Gamma \left( \alpha \right) }%
\int_{0}^{\infty }\frac{x^{\alpha -1}}{z^{-1}e^{x}-1}dx=\frac{z}{\Gamma
\left( \alpha \right) }\int_{0}^{1}\frac{\left( -\log u\right) ^{\alpha -1}}{%
1-uz}du
\end{equation*}
showing, by Hausdorff representation, that 
\begin{equation*}
\mathbf{P}\left( X_{\infty }>i\right) =\left( i+1\right) ^{-\alpha
}=\int_{0}^{1}u^{i}\lambda \left( du\right) \text{ where }\lambda \left(
du\right) =\frac{1}{\Gamma \left( \alpha \right) }\left( -\log u\right)
^{\alpha -1}du\text{ }
\end{equation*}
is the probability density of $U=e^{-X}$, with $X\overset{d}{\sim }$Gamma$%
\left( \alpha ,1\right) .$ The law of $X_{\infty }\geq 1$ is completely
monotone (and $X_{\infty }-1$ is ID). Note 
\begin{eqnarray*}
\Phi _{\infty }\left( z\right) &=&1-z^{-1}\left( 1-z\right) L_{\alpha
}\left( z\right) =z\int_{0}^{1}\frac{1-u}{1-uz}\lambda \left( du\right)
=z\Psi _{\infty }\left( z\right) \\
\Psi _{\infty }\left( z\right) &=&\mathbf{E}\left( z^{X_{\infty }-1}\right)
=\int_{0}^{1}\frac{1}{1-zu}\mu \left( du\right) \text{ where }\mu \left(
du\right) =\left( 1-u\right) \lambda \left( du\right)
\end{eqnarray*}
\newline

- \textbf{Zipf }($\alpha >1$): Suppose $\mathbf{P}\left( X_{\infty
}=i\right) =i^{-\alpha }/\varsigma \left( \alpha \right) $ with associated
pgf $\Phi _{\infty }\left( z\right) =L_{\alpha }\left( z\right) /L_{\alpha
}\left( 1\right) ,$ $L_{\alpha }\left( 1\right) =\varsigma \left( \alpha
\right) .$ Lagrange inversion formula gives the power-series expansion of $%
\Phi _{\infty }^{-1}\left( z\right) $. We have 
\begin{equation*}
\mathbf{P}\left( \nu \leq j\right) =\Phi _{\infty }^{-1}\left( 1-\mathbf{P}%
\left( X_{\infty }>i\right) \right)
\end{equation*}
where, with $\lambda _{0}\left( du\right) =\frac{1}{\Gamma \left( \alpha
\right) }\left( -\log u\right) ^{\alpha -1}du$%
\begin{eqnarray*}
\mathbf{P}\left( X_{\infty }>i\right) &=&\frac{1}{\varsigma \left( \alpha
\right) }\sum_{j>i}j^{-\alpha }=\frac{1}{\varsigma \left( \alpha \right) }%
\int_{0}^{1}\sum_{j>i}u^{j-1}\lambda _{0}\left( du\right) =\frac{1}{%
\varsigma \left( \alpha \right) }\int_{0}^{1}u^{i}\left( 1-u\right)
^{-1}\lambda _{0}\left( du\right) \\
&=&\int_{0}^{1}u^{i}\lambda \left( du\right) \text{ where }\lambda \left(
du\right) =\frac{\left( 1-u\right) ^{-1}}{\varsigma \left( \alpha \right)
\Gamma \left( \alpha \right) }\left( -\log u\right) ^{\alpha -1}du\text{ }
\end{eqnarray*}
is the probability density of $U=e^{-X}$, with $X$ having density 
\begin{equation*}
\frac{1}{\varsigma \left( \alpha \right) \Gamma \left( \alpha \right) }\frac{%
e^{-x}x^{\alpha -1}}{1-e^{-x}}\text{, }x>0.
\end{equation*}
So $X_{\infty }$ (and $X_{\infty }-1$) is CM. Thus $X_{\infty }-1$ is
infinitely divisible and even self-decomposable, say SD (see Example $12.18$
page $435$ of \cite{SH}). \newline

- \textbf{The critical case when} $X_{\infty }$ \textbf{has no moments of
any positive order}:

\begin{proposition}
Suppose that with $\beta >0$ and $L_{1}\left( x\right) =\log \left(
1+x\right) >0$, slowly varying at $\infty $ 
\begin{equation*}
\mathbf{P}\left( X_{\infty }=j\right) =\frac{C_{0}}{jL_{1}\left( j\right)
^{\beta +1}},\text{ }j\geq 1
\end{equation*}
where $C_{0}>0$ is the normalizing constant. Then $\mathbf{E}\left(
X_{\infty }^{q}\right) =\infty $ for all $q>0.$ In this case, $\mathbf{P}%
\left( X_{\infty }>j\right) \sim C_{0}\cdot L_{1}\left( j\right) ^{-\beta }$
as $j\rightarrow \infty $ with tails heavier than any power-law. Then:

$\left( i\right) $ The rv $\nu $ whose distribution solves (\ref{FE}) (as
from Proposition $2$) is a well-defined rv obeying $j\mathbf{P}\left( \nu
>j\right) \rightarrow e^{-\gamma }$ as $j\rightarrow \infty .$

$\left( ii\right) $ Furthermore 
\begin{equation*}
\mathbf{P}\left( \nu >j\right) \underset{j\uparrow \infty }{\sim }e^{-\gamma
}/j+d/\left( j\log j\right) +o\left( 1/\left( j\log j\right) \right) 
\end{equation*}
with 
\begin{equation*}
d=-\frac{\left( \beta +1\right) e^{-\gamma }\pi ^{2}}{12}<-\frac{e^{-\gamma
}\pi ^{2}}{12}.
\end{equation*}
By (\ref{L3}), the corresponding Lamperti chain is critical but it remains
positive recurrent for all $\beta >0.$
\end{proposition}

\emph{Proof:} $\left( i\right) $ This model for $X_{\infty }$ is indeed
obtained in the limit $\alpha \rightarrow 0$ of the ansatz ($\alpha >0$) 
\begin{equation*}
\mathbf{P}\left( X_{\infty }=j\right) =\frac{C_{0}}{j^{\alpha +1}L_{1}\left(
j\right) ^{\beta +1}},\text{ }i\geq 1,
\end{equation*}
extending the previous Sibuya example with tail index $\alpha $.

$\left( ii\right) $ In such an example of $X_{\infty }$ with logarithmic
tails, we have more precisely 
\begin{equation*}
\Phi _{\infty }\left( z\right) \underset{z\uparrow 1}{\sim }1-\frac{C_{0}}{%
\left( -\log \left( 1-z\right) \right) ^{\beta }}
\end{equation*}
with local inverse: $\Phi _{\infty }^{-1}\left( z\right) \underset{z\uparrow
1}{\sim }1-e^{-\left( \frac{1-z}{C_{0}}\right) ^{-1/\beta }}$. We get 
\begin{equation*}
\frac{1-\Phi _{\infty }\left( z\right) }{1-z}\underset{z\uparrow 1}{\sim }%
\frac{1}{1-z}\frac{C_{0}}{\left( -\log \left( 1-z\right) \right) ^{\beta }}
\end{equation*}
so that \cite{Fla}, with $C_{k}=\left( \frac{1}{\Gamma \left( \alpha \right) 
}\right) ^{\left( k\right) }\mid _{\alpha =1}$(in particular $C_{1}=\gamma $%
, $C_{2}=\gamma ^{2}-\pi ^{2}/6,$ with $C_{1}^{2}-C_{2}=\pi ^{2}/6$) 
\begin{equation*}
\mathbf{P}\left( X_{\infty }>j\right) \underset{j\uparrow \infty }{\sim }%
\frac{C_{0}}{\log ^{\beta }j}\left( 1-\frac{\beta C_{1}}{\log j}+\frac{\beta
\left( \beta +1\right) C_{2}}{2\log ^{2}j}+o\left( \frac{1}{\log ^{2}j}%
\right) \right) .
\end{equation*}
Observing $\left( 1-\frac{\beta C_{1}}{\log j}+\frac{\beta \left( \beta
+1\right) C_{2}}{2\log ^{2}j}\right) ^{-1/\beta }\underset{j\uparrow \infty 
}{\sim }1+\frac{C_{1}}{\log j}+\frac{\left( \beta +1\right) }{2\left( \log
j\right) ^{2}}\left( C_{1}^{2}-C_{2}\right) $, we are led to 
\begin{eqnarray*}
&&\mathbf{P}\left( \nu \leq j\right) \underset{j\uparrow \infty }{\sim }\Phi
_{\infty }^{-1}\left( 1-\mathbf{P}\left( X_{\infty }>j\right) \right) 
\underset{j\uparrow \infty }{\sim }1-\left( \frac{1}{j}\right) ^{\left( 1-%
\frac{\beta C_{1}}{\log j}+\frac{\beta \left( \beta +1\right) C_{2}}{2\log
^{2}j}\right) ^{-1/\beta }} \\
&&\underset{j\uparrow \infty }{\sim }1-e^{-\gamma }/j-d/\left( j\log
j\right) +o\left( 1/\left( j\log j\right) \right)
\end{eqnarray*}
with 
\begin{equation*}
d=-\frac{\left( \beta +1\right) e^{-\gamma }\pi ^{2}}{12}.
\end{equation*}
Because $d<-\pi ^{2}e^{-\gamma }/12$ for all $\beta >0,$ we conclude that $%
\left\{ X_{n}\right\} $ generated by this $\nu $ just remains always
positive-recurrent. $\Box $\newline

\textbf{- Null-recurrent issues.}\emph{\ }

Irreducible aperiodic Markov chains may have or not a non-trivial invariant
positive (infinite) measure, \cite{Harr}.

\begin{proposition}
In the null-recurrent case from (\ref{L2}), the Lamperti model has a non
trivial ($\neq \mathbf{0}$) invariant positive measure.
\end{proposition}

\emph{Proof:} To see a transition positive/null recurrence transition in the
critical case, suppose $\delta \left( j\right) :=\Delta F_{\infty }\left(
j\right) >0$ with $\Delta F_{\infty }\left( j\right) \rightarrow 0$ as $%
j\rightarrow \infty ,$ $\Phi _{\infty }\left( z\right) =\sum_{j\geq 1}\Delta
F_{\infty }\left( j\right) z^{j}$ convergent for all $z\in \left[ 0,1\right) 
$, $\Phi _{\infty }\left( 0\right) =0$ and $\Phi _{\infty }\left( 1\right)
=\infty .$ In this case $\Delta F_{\infty }\left( j\right) $ no longer is a
probability mass at $j$. One can search solutions of (\ref{FE}) in this case
as well and Proposition $2$ applies simply while substituting $\delta \left(
j\right) $ to $\pi \left( j\right) $ in the obtained expression of $\mathbf{P%
}\left( \nu \leq j\right) .$ Because $\mathbf{P}\left( \nu \leq j\right) $
only depends on the ratio $F_{\infty }\left( j\right) /F_{\infty }\left(
1\right) $, regardless of any normalization, such a sequence $\delta \left(
j\right) $ defines an invariant positive and infinite measure in the
null-recurrent case. $\Box $

The simplest example is the following: $\Delta F_{\infty }\left( j\right)
=1/j$, with $\Phi _{\infty }\left( z\right) =\sum_{j\geq 1}\Delta F_{\infty
}\left( j\right) z^{j}$ obeying $\Phi _{\infty }\left( 1\right) =\infty .$
We have $\Phi _{\infty }\left( z\right) =-\log \left( 1-z\right) $ so that,
upon inverting $\Phi _{\infty }$%
\begin{equation*}
\mathbf{P}\left( \nu \leq j\right) =1-e^{-\sum_{k=1}^{j}\frac{1}{k}}
\end{equation*}
a true pdf. Recalling $\sum_{k=1}^{j}\frac{1}{k}-\gamma -\log j\sim 1/\left(
2j\right) $, we get $\mathbf{P}\left( \nu >j\right) \sim e^{-\gamma
}/j+O\left( j^{-2}\right) .$ The constant $d$ in (\ref{L3}) is $d=0$ and the
Lamperti chain with a branching number $\nu $ distributed as such is
null-recurrent. This is also true if $\Delta F_{\infty }\left( j\right)
=1/\left[ j\log \left( 1+j\right) ^{\beta +1}\right] $ with $\beta <0$ or $%
\Delta F_{\infty }\left( j\right) =j^{-\alpha },$ $\alpha \in \left(
0,1\right) ,$ both expressions leading to a diverging series $\Phi _{\infty
}\left( 1\right) $.\newline

- \textbf{Transient issues: non-unicity of the invariant measure. }Whenever%
\textbf{\ }$\left\{ X_{n}\right\} $ is transient, one obvious solution to
the invariant measure equation $\mathbf{\pi }^{\prime }=\mathbf{\pi }%
^{\prime }P$ is $\mathbf{\pi }=\mathbf{0}$. This corresponds to the fact
that $X_{\infty }\overset{d}{\sim }\delta _{\infty }$. However this solution
is not unique and there are other invariant positive measures. The question
of the existence of a non-trivial invariant measure for transient chains was
raised by Harris, \cite{Har}.

To exhibit such an invariant measure, suppose $\delta \left( j\right)
:=\Delta F_{\infty }\left( j\right) >0$ with $\Phi _{\infty }\left( z\right)
=\sum_{j\geq 1}\Delta F_{\infty }\left( j\right) z^{j}$ convergent for all $%
z\in \left[ 0,1\right) $, $\Phi _{\infty }\left( 0\right) =0$ and $\Phi
_{\infty }\left( 1\right) =\infty .$ In this case $\Delta F_{\infty }\left(
j\right) $ no longer is a probability mass at $j$ either but it is no longer
required $\Delta F_{\infty }\left( j\right) \rightarrow 0$ as $j\rightarrow
\infty .$

\begin{proposition}
In the transient case from (\ref{L3}), the Lamperti model has a non trivial (%
$\neq \mathbf{0}$) invariant positive measure.
\end{proposition}

\emph{Proof:} One can search solutions of (\ref{FE}) in this case as well
and Proposition $2$ applies simply while substituting $\delta \left(
j\right) $ to $\pi \left( j\right) $ in the obtained expression of $\mathbf{P%
}\left( \nu \leq j\right) .$ Because, from (\ref{FES}), $\mathbf{P}\left(
\nu \leq j\right) $ only depends on the ratio $F_{\infty }\left( j\right)
/F_{\infty }\left( 1\right) $ regardless of any normalization, such a
sequence $\delta \left( j\right) $ defines an invariant measure in the
transient case as well. $\Box $

- The simplest explicit example is the following counting measure one: $%
\delta \left( j\right) =\Delta F_{\infty }\left( j\right) =1$, $F_{\infty
}\left( j\right) =j$, with $\Phi _{\infty }\left( z\right) =\sum_{j\geq
1}\Delta F_{\infty }\left( j\right) z^{j}=z/\left( 1-z\right) $ obeying $%
\Phi _{\infty }\left( 1\right) =\infty .$ There is a solution to (\ref{FE})
which is 
\begin{equation*}
\mathbf{P}\left( \nu \leq j\right) =\Phi _{\infty }^{-1}\left( j\right) =%
\frac{j}{1+j}.
\end{equation*}
We have: $\mathbf{P}\left( \nu >j\right) =1/\left( 1+j\right) $ so that $j%
\mathbf{P}\left( \nu >j\right) \underset{j\rightarrow \infty }{\rightarrow }%
1>e^{-\gamma }$, indeed corresponding to a transient case.

- Suppose now $\Delta F_{\infty }\left( j\right) =j$, $F_{\infty }\left(
j\right) =j\left( j+1\right) /2,$ so with $\Phi _{\infty }\left( z\right)
=\sum_{j\geq 1}\Delta F_{\infty }\left( j\right) z^{j}=z/\left( 1-z\right)
^{2}$ obeying $\Phi _{\infty }\left( 1\right) =\infty .$ There is a solution
to (\ref{FE}) which is 
\begin{equation*}
\mathbf{P}\left( \nu \leq j\right) =\Phi _{\infty }^{-1}\left( \frac{j\left(
j+1\right) }{2}\right) =\frac{j\left( j+1\right) +1-\sqrt{1+2j\left(
j+1\right) }}{j\left( j+1\right) }.
\end{equation*}
When inverting $\Phi _{\infty }\left( z\right) $ we have chosen the branch
for which $\Phi _{\infty }^{-1}\left( 0\right) =0$. We have: $\mathbf{P}%
\left( \nu >j\right) =\left( \sqrt{1+2j\left( j+1\right) }-1\right) /\left(
j\left( j+1\right) \right) $ so that $j\mathbf{P}\left( \nu >j\right)
\rightarrow \sqrt{2}>e^{-\gamma }$, also corresponding to a transient case.
Defining the reversed failure rate of the sequence $\delta \left( j\right) $
as 
\begin{equation*}
\overline{r}\left( j\right) =\frac{\delta \left( j\right) }{%
\sum_{k=1}^{j}\delta \left( k\right) }=\frac{\Delta F_{\infty }\left(
j\right) }{F_{\infty }\left( j\right) }\text{, }j\geq 1,
\end{equation*}
we conclude that in both examples, $\overline{r}\left( j\right) \asymp 1/j$
so with decreasing reversed failure rate.\newline

\emph{Remark:} By the ergodic theorem:

- in case (\ref{L1}): 
\begin{equation*}
n^{-1}\sum_{m=1}^{n}\mathbf{1}\left( X_{m}=j\mid X_{0}\overset{d}{\sim }%
\mathbf{\pi }_{0}\right) \rightarrow \pi \left( j\right) \text{ as }%
n\rightarrow \infty ,
\end{equation*}

- in cases (\ref{L2}) and (\ref{L3}): For all states $i,j\geq 1$%
\begin{equation*}
\frac{\sum_{m=1}^{n}\mathbf{1}\left( X_{m}=i\mid X_{0}\overset{d}{\sim }%
\mathbf{\pi }_{0}\right) }{\sum_{m=1}^{n}\mathbf{1}\left( X_{m}=j\mid X_{0}%
\overset{d}{\sim }\mathbf{\pi }_{0}\right) }\rightarrow \frac{\delta \left(
i\right) }{\delta \left( j\right) }\text{ as }n\rightarrow \infty .\text{ }%
\diamondsuit 
\end{equation*}
\newline

- \textbf{Poisson target: }We finally develop some additional examples in
the recurrent case, not in the latter classes and related to the fundamental
Poisson distribution class:

- \textbf{Shifted} \textbf{Poisson: }

Suppose\textbf{\ }$\Phi _{\infty }\left( z\right) =\mathbf{E}\left(
z^{X_{\infty }}\right) =ze^{\lambda \left( z-1\right) }=z\Psi _{\infty
}\left( z\right) ,$ \textbf{(}$\Psi _{\infty }\left( z\right) $ is the pgf
of an ID Poisson rv which is log-concave). Then 
\begin{eqnarray*}
\frac{1}{n}\left[ z^{n-1}\right] \Psi _{\infty }\left( z\right) ^{-n} &=&%
\frac{1}{n}\left[ z^{n-1}\right] e^{-n\lambda \left( z-1\right) }=\left(
-1\right) ^{n-1}\frac{e^{\lambda n}}{n!}\left( n\lambda \right) ^{n-1} \\
F_{\infty }\left( j\right) &=&e^{-\lambda }\sum_{k=0}^{j-1}\frac{\lambda ^{k}%
}{k!} \\
F\left( j\right) &=&\sum_{n\geq 1}\left( -1\right) ^{n-1}\frac{e^{\lambda n}%
}{n!}\left( n\lambda \right) ^{n-1}F_{\infty }\left( j\right) ^{n} \\
&=&\sum_{n\geq 1}\left( -1\right) ^{n-1}\frac{\left( n\lambda \right) ^{n-1}%
}{n!}\left( \sum_{k=0}^{j-1}\frac{\lambda ^{k}}{k!}\right) ^{n}=W_{\lambda
}\left( \sum_{k=0}^{j-1}\frac{\lambda ^{k}}{k!}\right)
\end{eqnarray*}
The Lambert function, solving $x=W\left( x\right) e^{W\left( x\right) },$ is
(by Lagrange inversion formula): 
\begin{eqnarray*}
W\left( x\right) &=&\sum_{n\geq 1}\left( -1\right) ^{n-1}\frac{n^{n-1}}{n!}%
x^{n}\text{ hence} \\
W_{\lambda }\left( x\right) &:&=\lambda ^{-1}W\left( \lambda x\right)
=\sum_{n\geq 1}\left( -1\right) ^{n-1}\frac{\left( n\lambda \right) ^{n-1}}{%
n!}x^{n}.
\end{eqnarray*}
And $W_{\lambda }\left( x\right) $ solves: $x=W_{\lambda }\left( x\right)
e^{\lambda W_{\lambda }\left( x\right) }$. It is positive and increasing
when $x>0$, so $F\left( j\right) $ is a well-defined pdf if $F\left( \infty
\right) =W_{\lambda }\left( e^{\lambda }\right) =1,$ which is the case.%
\newline

- \textbf{Poisson conditioned to be positive: }

Suppose\textbf{\ }$\Phi _{\infty }\left( z\right) =\mathbf{E}\left(
z^{X_{\infty }}\right) =\left( e^{\lambda z}-1\right) /\left( e^{\lambda
}-1\right) $, leading directly to $\Phi _{\infty }^{-1}\left( z\right) =%
\frac{1}{\lambda }\log \left( 1+z\left( e^{\lambda }-1\right) \right) .$
Then 
\begin{eqnarray*}
F_{\infty }\left( j\right) &=&\frac{1}{e^{\lambda }-1}\sum_{k=1}^{j}\frac{%
\lambda ^{k}}{k!} \\
F\left( j\right) &=&\Phi _{\infty }^{-1}\left( F_{\infty }\left( j\right)
\right) =\frac{1}{\lambda }\log \left( 1+F_{\infty }\left( j\right) \left(
e^{\lambda }-1\right) \right) ,
\end{eqnarray*}
which defines a pdf with $F\left( \infty \right) =1$. In this case, although 
$\Psi _{\infty }\left( z\right) =z^{-1}\Phi _{\infty }\left( z\right) $ is
not the pgf of an ID rv, the calculation of $F\left( j\right) $ is
straightforward.

\subsection{Examples of $\nu \rightarrow \nu _{\left( N\right) }$ with
finite support $\left\{ 1,...,N\right\} $.}

In this Sub-section, we look at situations where both $\left( X_{\infty
},\nu _{\left( N\right) }\right) $ have finite support $\left\{
1,...,N\right\} $. Note that if $\nu $ has support $\left\{ 1,...,N\right\} $%
, so does $\left\{ X_{n}\right\} $ (defined recursively by $%
X_{n+1}=\max_{j=1,...,X_{n}}\nu _{j,n+1}$) and then $X_{\infty }$.
Conversely, if $X_{\infty }$ has support $\left\{ 1,...,N\right\} $, there
exists $\nu $ with support $\left\{ 1,...,N\right\} $ such that $%
X_{n+1}=\max_{j=1,...,X_{n}}\nu _{j,n+1}$ defines a sequence $\left(
X_{n}\right) $ with finite support. In such cases, the Lamperti Markov chain
will always be ergodic in view of its transition matrix $P_{\left( N\right)
} $ being irreducible. We shall let $\pi _{\left( N\right) }\left( k\right) =%
\mathbf{P}\left( X_{\infty }=k\right) .$

- \textbf{The general case:}

Suppose $\Phi _{\infty }\left( z\right) =\sum_{k=1}^{N}\pi _{\left( N\right)
}\left( k\right) z^{k}$, so that $\Psi _{\infty }\left( z\right)
=\sum_{k=0}^{N-1}\pi _{\left( N\right) }\left( k+1\right) z^{k}$. We have 
\begin{eqnarray*}
\Psi _{\infty }\left( z\right) ^{-\alpha } &=&\pi _{\left( N\right) }\left(
1\right) ^{-\alpha }\left( 1+\sum_{k=1}^{N-1}\frac{\pi _{\left( N\right)
}\left( k+1\right) }{\pi _{\left( N\right) }\left( 1\right) }z^{k}\right)
^{-\alpha } \\
&=&\pi _{\left( N\right) }\left( 1\right) ^{-\alpha }\sum_{l\geq
0}z^{l}\sum_{k=0}^{l}\left( -1\right) ^{k}\left[ \alpha \right]
_{k}\sum_{{}}^{*}\prod_{m=1}^{N-1}\frac{\left( \pi _{\left( N\right) }\left(
m+1\right) /\pi _{\left( N\right) }\left( 1\right) \right) ^{k_{m}}}{k_{m}!}
\end{eqnarray*}
where the star sum runs over $k_{m}\geq 0$, $m=1,...,N-1$ obeying $%
\sum_{m=1}^{N-1}k_{m}=k$ and $\sum_{m=1}^{N-1}mk_{m}=l$. From this, we
obtain the finite support version of (\ref{FES}) as

\begin{proposition}
For any given $X_{\infty }$ with support $\left\{ 1,...,N\right\} $, the
mapping $X_{\infty }\rightarrow \nu _{\left( N\right) }$ is one-to-one and
onto. With 
\begin{equation*}
C_{n-1,0}^{\left( N-1\right) }=\delta _{n-1,0}\text{ and }C_{n-1,k}^{\left(
N-1\right) }:=\sum_{{}}^{*}\prod_{m=1}^{N-1}\frac{\left( \pi _{\left(
N\right) }\left( m+1\right) /\pi _{\left( N\right) }\left( 1\right) \right)
^{k_{m}}}{k_{m}!}
\end{equation*}
where the star sum runs over $k_{m}\geq 0$, $m=1,...,N-1$ obeying $%
\sum_{m=1}^{N-1}k_{m}=k$ and $\sum_{m=1}^{N-1}mk_{m}=n-1\geq k$, 
\begin{equation}
\varphi _{n}:=\left[ z^{n}\right] \Phi _{\infty }^{-1}\left( z\right) =\frac{%
1}{n}\left[ z^{n-1}\right] \Psi _{\infty }\left( z\right) ^{-n}=\frac{\pi
_{\left( N\right) }\left( 1\right) ^{-n}}{n}\sum_{k=0}^{n-1}\left( -1\right)
^{k}\left[ n\right] _{k}C_{n-1,k}^{\left( N-1\right) }.  \label{FS1}
\end{equation}
So, with $h_{1}=1$ and $h_{n}=\frac{1}{n}\sum_{k=1}^{n-1}\left( -1\right)
^{k}\left[ n\right] _{k}C_{n-1,k}^{\left( N-1\right) },$ $n\geq 2$%
\begin{equation}
\mathbf{P}\left( \nu _{\left( N\right) }\leq j\right) =\sum_{n\geq
1}h_{n}\cdot \left( \mathbf{P}\left( X_{\infty }\leq j\right) /\pi _{\left(
N\right) }\left( 1\right) \right) ^{n}  \label{FS2}
\end{equation}
is the pdf of $\nu _{\left( N\right) }$ associated to any $\mathbf{P}\left(
X_{\infty }\leq j\right) =\sum_{k=1}^{j}\pi _{\left( N\right) }\left(
k\right) $, $j=1,...,N$ obeying $\mathbf{P}\left( X_{\infty }\leq N\right)
=1.$
\end{proposition}

\emph{Remark: }

$\left( i\right) $ $\varphi _{1}=1/\pi _{\left( N\right) }\left( 1\right) $ (%
$h_{1}=1$) and for $n\geq 2$, the sum over $k$ giving the expression of $%
\varphi _{n}$ (or of $h_{n}$) can start at $k=1$.

$\left( ii\right) $ if (a separable case in $\left( k,N\right) $): $\pi
_{\left( N\right) }\left( k\right) =a_{k}/A_{N},$ $a_{k}\geq 0$, where $%
A_{N}=\sum_{k=1}^{N}a_{k}$ is a normalization factor, the law of $\nu
_{\left( N\right) }$ does not depend on $A_{N}$ because it only depends on
the ratios $\pi _{\left( N\right) }\left( k\right) /\pi _{\left( N\right)
}\left( 1\right) =a_{k}/a_{1}.$ $\diamondsuit $ \newline

\textbf{Examples: }Just like in the infinite-dimensional case, there are
examples amenable to a straightforward calculation.

$\left( i\right) $ Suppose

\begin{equation*}
\Phi _{\infty }\left( z\right) =\frac{\left( q+pz\right) ^{N}-q^{N}}{1-q^{N}}%
,
\end{equation*}
corresponding to a binomial model restricted to $\left\{ 1,...,N\right\} $
with 
\begin{eqnarray*}
\mathbf{P}\left( X_{\infty }=k\right) &=&\left[ z^{k}\right] \Phi _{\infty
}\left( z\right) =\frac{1}{1-q^{N}}\binom{N}{k}p^{k}q^{N-k} \\
\mathbf{P}\left( X_{\infty }\leq j\right) &=&\sum_{k=1}^{j}\mathbf{P}\left(
X_{\infty }=k\right) .
\end{eqnarray*}
By direct inversion of $\Phi _{\infty }\left( z\right) $, we have that 
\begin{equation*}
\Phi _{\infty }^{-1}\left( \mathbf{P}\left( X_{\infty }\leq j\right) \right)
=\frac{q}{p}\left[ \left( 1+\sum_{k=1}^{j}\binom{N}{k}\left( \frac{p}{q}%
\right) ^{k}\right) ^{1/N}-1\right]
\end{equation*}
is the pdf $\mathbf{P}\left( \nu _{\left( N\right) }\leq j\right) $ of some
rv $\nu _{\left( N\right) }.$ Note $\mathbf{P}\left( \nu _{\left( N\right)
}\leq N\right) =\Phi _{\infty }^{-1}\left( \mathbf{P}\left( X_{\infty }\leq
N\right) \right) =1.$\newline

$\left( ii\right) $ Suppose

\begin{equation*}
\Phi _{\infty }\left( z\right) =z\left( q+pz\right) ^{N-1}=z\Psi _{\infty
}\left( z\right)
\end{equation*}
corresponding to a shifted binomial model supported $\left\{ 1,...,N\right\} 
$%
\begin{eqnarray*}
\mathbf{P}\left( X_{\infty }=k\right) &=&\left[ z^{k}\right] \Phi _{\infty
}\left( z\right) =\binom{N-1}{k-1}p^{k-1}q^{N-k} \\
\mathbf{P}\left( X_{\infty }\leq j\right) &=&\sum_{k=1}^{j}\mathbf{P}\left(
X_{\infty }=k\right) .
\end{eqnarray*}
With $n\geq 1$, we have that $\Phi _{\infty }^{-1}\left( z\right)
=\sum_{n\geq 1}\frac{z^{n}}{n}\left[ z^{n-1}\right] \Psi _{\infty }\left(
z\right) ^{-n}$ with 
\begin{eqnarray*}
\varphi _{n} &=&\frac{q^{-n\left( N-1\right) }}{n}\left[ z^{n-1}\right]
\left( 1+\frac{p}{q}z\right) ^{-n\left( N-1\right) } \\
&=&\left( -1\right) ^{n-1}q^{-n\left( N-1\right) }\left( \frac{p}{q}\right)
^{n-1}\frac{\left[ n\left( N-1\right) \right] _{n-1}}{n!}
\end{eqnarray*}
\begin{equation*}
\mathbf{P}\left( \nu _{\left( N\right) }\leq j\right) =\sum_{n\geq 1}\varphi
_{n}\mathbf{P}\left( X_{\infty }\leq j\right) ^{n}
\end{equation*}
The rv $\nu _{\left( N\right) }$ has support $\left\{ 1,...,N\right\} $.%
\newline

$\left( iii\right) $ Truncation of the infinite-dimensional model.

This situation occurs if, for $\left( \pi _{\left( N\right) }\left( k\right)
,k=1,...,N\right) $, we consider the normalized restriction of the invariant
measure $\mathbf{\pi }$ with full support $\Bbb{N}$ to its $N$ first
entries. For example, assuming $\left( \pi \left( k\right) =pq^{k-1},k\geq
1\right) $ is geometric, we get 
\begin{equation*}
\mathbf{P}\left( \nu _{\left( N\right) }\leq j\right) =\sum_{n\geq
1}h_{n}\cdot \left( \frac{1-q^{j}}{p}\right) ^{n}
\end{equation*}
where $h_{n}=\frac{q^{n-1}}{n}\sum_{k=1}^{n-1}\left( -1\right) ^{k}\frac{%
\left[ n\right] _{k}}{k!}\sum_{{}}^{*}\frac{k!}{\prod_{m=1}^{N-1}k_{m}!}=%
\frac{q^{n-1}}{n}\sum_{k=1}^{n-1}\left( -1\right) ^{k}\frac{\left[ n\right]
_{k}}{k!}\left( N-1\right) ^{k}$, so that with $A_{n,N}=\sum_{k=1}^{n-1}%
\left( -1\right) ^{k}\frac{\left[ n\right] _{k}}{k!}\left( N-1\right) ^{k}$, 
\begin{equation*}
\mathbf{P}\left( \nu _{\left( N\right) }\leq j\right) =\frac{1}{q}%
\sum_{n\geq 1}\frac{A_{n,N}}{n}\cdot \left( \frac{q\left( 1-q^{j}\right) }{p}%
\right) ^{n}.\text{ }\diamondsuit
\end{equation*}

\begin{proposition}
Take any probability measure $\mathbf{\pi }_{\left( N\right) }$ with support 
$\left\{ 1,...,N\right\} $. Compute $F_{\left( N\right) }\left( j\right) =%
\mathbf{P}\left( \nu _{\left( N\right) }\leq j\right) $ from $\mathbf{\pi }%
_{\left( N\right) }$ as from (\ref{FS2}). Construct the $N\times N$
stochastic matrix $P_{\left( N\right) }$ with entries $P_{\left( N\right)
}\left( i,j\right) =F_{\left( N\right) }\left( j\right) ^{i}-F_{\left(
N\right) }\left( j-1\right) ^{i}$, $i,j\in \left\{ 1,...,N\right\} .$ The
matrix $P_{\left( N\right) }$ is the transition matrix of some ergodic
Lamperti chain $X_{n}^{\left( N\right) }$ with state-space $\left\{
1,...,N\right\} ,$ having $\mathbf{\pi }_{\left( N\right) }$ as invariant
probability measure and reproduction mechanism $\nu _{\left( N\right) }$.
The MC $\left\{ X_{n}^{\left( N\right) }\right\} $ is failure rate monotone.
Furthermore: 
\begin{equation*}
\mathbf{P}_{\mathbf{\pi }_{0}}\left( X_{n}^{\left( N\right) }=j\right) 
\underset{n,N\rightarrow \infty }{\rightarrow }\pi \left( j\right) 
\end{equation*}
\end{proposition}

\emph{Proof:} The reasons are similar to the ones raised for the Lamperti
chain taking values in $\Bbb{N}$. The probability $\mathbf{P}\left(
X_{n+1}^{\left( N\right) }\leq j\mid X_{n}^{\left( N\right) }=i\right)
=F_{\left( N\right) }\left( j\right) ^{i}$ is a decreasing function of $i$,
for all $j$: the MC $\left\{ X_{n}^{\left( N\right) }\right\} $ is
stochastically monotone. The cumulated transition matrix : $P_{\left(
N\right) }^{c}\left( i,j\right) =\sum_{k=1}^{j}P_{\left( N\right) }\left(
i,k\right) $ obeys:

\begin{equation*}
P_{\left( N\right) }^{c}\left( i_{1},j_{1}\right) P_{\left( N\right)
}^{c}\left( i_{2},j_{2}\right) \geq P_{\left( N\right) }^{c}\left(
i_{1},j_{2}\right) P_{\left( N\right) }^{c}\left( i_{2},j_{1}\right) ,
\end{equation*}
for all $i_{1}<i_{2}$ and $j_{1}<j_{2}$ ($P_{\left( N\right) }^{c}$ is
totally positive of order $2$): the MC $\left\{ X_{n}^{\left( N\right)
}\right\} $ is failure rate monotone.

Note the induced Kirchoff determinantal identities for finite matrices: $\pi
_{\left( N\right) }\left( j\right) =\det \left[ \left( I-P_{\left( N\right)
}\right) ^{\left( j,j\right) }\right] .$ The last statement is obvious. $%
\Box $

\begin{corollary}
(Truncation of $X_{n}$)

$\left( i\right) $ Take for $\mathbf{\pi }_{\left( N\right) }$ the
restriction to $\left\{ 1,...,N\right\} $\ of the invariant measure\emph{\ }$%
\mathbf{\pi }$ of the Lamperti model with countable state-space, so with: $%
\mathbf{\pi }_{\left( N\right) }\left( k\right) =\pi \left( k\right)
/\sum_{k=1,...,N}\pi \left( k\right) $\emph{, }$k=1,...,N$\emph{.}

$\left( ii\right) $ Take for $\mathbf{\pi }_{\left( N\right) }$ the
restriction to $\left\{ 1,...,N-1\right\} $\ of the invariant measure\emph{\ 
}$\mathbf{\pi }$ of the Lamperti model with countable state-space, so with: $%
\mathbf{\pi }_{\left( N\right) }\left( k\right) =\pi \left( k\right) $\emph{%
, }$k=1,...,N-1$, $\pi _{\left( N\right) }\left( N\right) =\sum_{k\geq N}\pi
\left( k\right) $\emph{.}

Constructing the corresponding transition matrices $P_{\left( N\right) }$,
in both cases, the truncations preserve the failure rate monotonicity of $P.$
\end{corollary}

The corresponding Lamperti chains $X_{n}^{\left( N\right) }$ with
state-space $\left\{ 1,...,N\right\} ,$ having $\mathbf{\pi }_{\left(
N\right) }$ as restricted invariant measure and reproduction mechanism $\nu
_{\left( N\right) }$ are called the truncated Lamperti chains up to state $N$%
.

\emph{Remarks:}

- The case $\left( i\right) $ is simpler because in this separable case, the
corresponding law of $\nu _{\left( N\right) }$ does not depend on the
normalization factor $\sum_{k=1,...,N}\pi \left( k\right) $.

- Censored Markov chain\textbf{\ }(\cite{ZL}, \cite{GH}): with $%
P_{11}=Q_{\left( N\right) }$ and 
\begin{equation*}
P=\left[ 
\begin{array}{ll}
P_{11} & P_{12} \\ 
P_{21} & P_{22}
\end{array}
\right] ,
\end{equation*}
define 
\begin{equation*}
P_{\left( N\right) }=P_{11}+P_{12}\left( I-P_{22}\right) ^{-1}P_{21}.
\end{equation*}
Let $Q_{2,2}=\left( I-P_{22}\right) ^{-1}$ be the fundamental matrix of $%
P_{22}$, with $Q_{2,2}\left( i,j\right) $ the mean number of visits to state 
$j$ in $\left\{ N+1,...,\infty \right\} $ starting from $i$ in $\left\{
N+1,...,\infty \right\} $, before visiting first $\left\{ 1,...,N\right\} $.
The matrix element $\left( P_{12}Q_{2,2}P_{21}\right) \left( i,j\right) $ is
the taboo probability of the paths from states $i$ to $j$ both in $\left\{
1,...,N\right\} $ which are not allowed to visit $\left\{ 1,...,N\right\} $
in between. $P_{\left( N\right) }$ has invariant measure $\mathbf{\pi }%
_{\left( N\right) }^{\prime }=\left( \pi _{1},...,\pi _{N}\right) /$norm
(the restriction of $\mathbf{\pi }$ to its $N$ first entries). However, it
is not clear that such a $P_{\left( N\right) }$\ is SM (probably not) nor
that $P_{\left( N\right) }^{c}$\ is FRM$.$\ Besides, $P_{\left( N\right) }$\
has a complicated structure in case of Lamperti. Truncating a Markov chain
invariant measure while preserving the monotonicity properties of the
original is not so straightforward.\emph{\ }$\diamondsuit $

\section{Brown's analysis of the truncated Lamperti model}

In this Section, we consider the truncated version $\left\{ X_{n}^{\left(
N\right) }\right\} $ of the chain $\left\{ X_{n}\right\} $ corresponding to
the one preserving the $N$ first entries of the full invariant measure $%
\mathbf{\pi }$ of $\left\{ X_{n}\right\} $, meaning $\pi \left( i\right)
\rightarrow \pi _{\left( N\right) }\left( i\right) =\pi \left( i\right)
/\sum_{i=1}^{N}\pi \left( i\right) $, $i=1,...,N$ (the restriction to $%
\left\{ 1,...,N\right\} $\ of the full invariant measure supported by $\Bbb{N%
}$). This MC has totally ordered state-space, with $\left\{ N\right\} $ as a
maximal element. It is a separable case and this truncation preserves the
failure-rate monotonicity of $P^{c}:$ $P_{\left( N\right) }^{c}$ remains
FRM, else $P_{\left( N\right) }^{c}$ is TP$_{2}$. As in \cite{Brown}, we
shall be concerned by the relationship existing between the first hitting
times of both state $\left\{ N\right\} $ and the restricted invariant
measure $\mathbf{\pi }_{\left( N\right) },$ given $X_{0}^{\left( N\right) }%
\overset{d}{\sim }\mathbf{\pi }_{0}$. We will assume $\pi _{0}\left(
N\right) =0$, to ensure that $\left\{ X_{n}^{\left( N\right) }\right\} $
hits $\left\{ N\right\} $ for the first time with positive probability after
at least one time unit. To illustrate his theory, Brown designs some ad hoc (%
$4\times 4$) FRM matrices; the truncated Lamperti chain is a more relevant
example. The following general results for hitting times hold for the
Lamperti truncated chain (see also \cite{Lorekth} for a survey).

\begin{proposition}
\cite{Brown}. Suppose $\mathbf{\pi }_{0}$ is such that $\pi _{0}\left(
i\right) /\pi _{\left( N\right) }\left( i\right) $ decreases with $i$ and $%
\pi _{0}\left( N\right) =0$. Then

$\left( i\right) $ $\mathbf{P}\left( X_{n}^{\left( N\right) }=N\mid
X_{0}^{\left( N\right) }\overset{d}{\sim }\mathbf{\pi }_{0}\right) $ is
non-decreasing with $n.$

$\left( ii\right) $ Let $\tau _{i,j}=\inf \left( n\geq 1:X_{n}^{\left(
N\right) }=j\mid X_{0}^{\left( N\right) }=i\right) $, with $\tau _{j,j}:=0$.
With $\tau _{\mathbf{\pi }_{0},j}=\inf \left( n\geq 1:X_{n}^{\left( N\right)
}=j\mid X_{0}^{\left( N\right) }\overset{d}{\sim }\mathbf{\pi }_{0}\right) :$%
\begin{equation}
\tau _{\mathbf{\pi }_{0},N}\overset{d}{=}T_{\left( N\right) }+\tau _{\mathbf{%
\pi }_{\left( N\right) },N}  \label{B0}
\end{equation}
where $T_{\left( N\right) }\geq 1$ and $\tau _{\mathbf{\pi }_{\left(
N\right) },N}\geq 0$ are independent.
\end{proposition}

\emph{Proof:} The condition that $\mathbf{\pi }_{0}$ is such that $\pi
_{0}\left( i\right) /\pi _{\left( N\right) }\left( i\right) $ decreases with 
$i$ holds if $\pi _{0}\left( i\right) =\delta _{i,1}$ and also if $\pi
_{0}\left( i\right) =z^{i}\pi _{\left( N\right) }\left( i\right) /$norm, $%
i=1,...,N-1$ for some $z\in \left( 0,1\right) $). It says that the initial
probability mass assigned to states near the bottom state $\left\{ 1\right\} 
$ should exceed the one assigned by $\mathbf{\pi }_{\left( N\right) }$. In
particular: $\pi _{0}\left( 1\right) >\pi _{\left( N\right) }\left( 1\right)
.$

The proof of this statement was derived in \cite{Brown} in a continuous-time
setting and is easily adaptable to discrete-time.

$\left( i\right) $ Let $\mathbf{e}_{N}^{\prime }=\left( 0,...,0,1\right) $
be an $N-$dimensional row vector, with $1$ in position $N$. Because $%
P_{\left( N\right) }^{c}$ is FRM (in particular SM), $\mathbf{P}\left(
X_{n}^{\left( N\right) }=N\mid X_{0}^{\left( N\right) }\overset{d}{\sim }%
\mathbf{\pi }_{0}\right) =\mathbf{\pi }_{0}^{\prime }P_{\left( N\right) }^{n}%
\mathbf{e}_{N}$ is non-decreasing with $n$ \cite{Brown}. See Lemma $4.2$ of 
\cite{Brown} where it is shown that this condition is fulfilled if $%
\overline{P}_{\left( N\right) }^{n}\left( \mathbf{\pi }_{0},j\right) /%
\overline{\pi }_{\left( N\right) }\left( j\right) $ is decreasing in $j$,
which is the case for FRM Markov chains. Here $\overline{\pi }_{\left(
N\right) }\left( j\right) =\sum_{k=j}^{N}\overline{\pi }_{\left( N\right)
}\left( k\right) $ and $\overline{P}_{\left( N\right) }^{n}\left( \mathbf{%
\pi }_{0},j\right) =\sum_{k=j}^{N}P_{\left( N\right) }^{n}\left( \mathbf{\pi 
}_{0},k\right) .$ It is needed in the proof that $\pi _{0}\left( i\right)
/\pi _{\left( N\right) }\left( i\right) $ decreases with $i.$

$\left( ii\right) $ Owing to $P_{\left( N\right) }^{n}\left( \mathbf{\pi }%
_{0},N\right) =\mathbf{\pi }_{0}^{\prime }P_{\left( N\right) }^{n}\mathbf{e}%
_{N}\rightarrow \pi _{\left( N\right) }\left( N\right) $ as $n\rightarrow
\infty $, $\mathbf{\pi }_{0}^{\prime }P_{\left( N\right) }^{n}\mathbf{e}%
_{N}/\pi _{\left( N\right) }\left( N\right) $ is a probability distribution
function of some rv $T_{\left( N\right) }$ with 
\begin{equation*}
\mathbf{P}\left( T_{\left( N\right) }\leq n\right) =\mathbf{\pi }%
_{0}^{\prime }P_{\left( N\right) }^{n}\mathbf{e}_{N}/\pi _{\left( N\right)
}\left( N\right) \text{, }n\geq 0.
\end{equation*}
With $G_{\mathbf{\pi }_{0},N}\left( z\right) =\mathbf{\pi }_{0}^{\prime
}\sum_{n\geq 0}z^{n}P_{\left( N\right) }^{n}\mathbf{e}_{N}=\mathbf{\pi }%
_{0}^{\prime }\left( I-zP_{\left( N\right) }\right) ^{-1}\mathbf{e}_{N}$, a
Green kernel of $P_{\left( N\right) }$, we thus have 
\begin{equation}
\mathbf{E}\left( z^{T_{\left( N\right) }}\right) =\frac{1-z}{\pi _{\left(
N\right) }\left( N\right) }G_{\mathbf{\pi }_{0},N}\left( z\right) .
\label{B1}
\end{equation}

Now 
\begin{equation*}
\mathbf{\pi }_{\left( N\right) }\left( N\right) =\mathbf{P}_{\mathbf{\pi }%
_{\left( N\right) }}\left( X_{n}^{\left( N\right) }=N\right)
=\sum_{m=0}^{n}P_{\left( N\right) }^{n-m}\left( N,N\right) \mathbf{P}\left(
\tau _{\mathbf{\pi }_{\left( N\right) },N}=m\right) ,
\end{equation*}
of convolution type. Taking the generating function of both sides 
\begin{equation}
\mathbf{E}\left( z^{\tau _{\mathbf{\pi }_{\left( N\right) },N}}\right) =%
\frac{\mathbf{\pi }_{\left( N\right) }\left( N\right) }{\left( 1-z\right)
G_{N,N}\left( z\right) },  \label{B2}
\end{equation}
where $G_{N,N}\left( z\right) =\sum_{m\geq 0}z^{m}P_{\left( N\right)
}^{m}\left( N,N\right) =\left( I-zP_{\left( N\right) }\right) ^{-1}\left(
N,N\right) $ is the Green kernel of $\left\{ X_{n}^{\left( N\right)
}\right\} $ at $\left( N,N\right) .$ Similarly, 
\begin{equation*}
\mathbf{P}_{\mathbf{\pi }_{0}}\left( X_{n}^{\left( N\right) }=N\right)
=\sum_{m=0}^{n}P_{\left( N\right) }^{n-m}\left( N,N\right) \mathbf{P}\left(
\tau _{\mathbf{\pi }_{0},N}=m\right)
\end{equation*}
leading to, 
\begin{equation}
G_{\mathbf{\pi }_{0},N}\left( z\right) =G_{N,N}\left( z\right) \mathbf{E}%
\left( z^{\tau _{\mathbf{\pi }_{0},N}}\right)  \label{B3}
\end{equation}
Taking the product of (\ref{B2}-\ref{B3}), and recalling (\ref{B1}), we get 
\begin{equation}
\phi _{\mathbf{\pi }_{0},N}\left( z\right) :=\mathbf{E}\left( z^{\tau _{%
\mathbf{\pi }_{0},N}}\right) =\mathbf{E}\left( z^{T_{\left( N\right)
}}\right) \mathbf{E}\left( z^{\tau _{\mathbf{\pi }_{\left( N\right)
},N}}\right) .\text{ }\Box  \label{B3a}
\end{equation}

The latter equation indicates that $\tau _{\mathbf{\pi }_{0},N}\geq 1$ is
stochastically larger than $\tau _{\mathbf{\pi }_{\left( N\right) },N}$: it
takes a shorter time for $\left\{ X_{n}^{\left( N\right) }\right\} $ to
first hit $\left\{ N\right\} $ starting from $\mathbf{\pi }_{\left( N\right)
}$ than starting from $\mathbf{\pi }_{0}.$ The time to first hit state $%
\left\{ N\right\} $ is important in the Lamperti context because at this
instant, the progeny after selection is the maximum possible. But of course
the process will not remain in that state unless one forces the chain to
have $\left\{ N\right\} $ absorbing.

As a result also, $T_{\left( N\right) }$ interprets as $\tau _{\mathbf{\pi }%
_{0},\mathbf{\pi }_{\left( N\right) }},$ the first hitting time of $\mathbf{%
\pi }_{\left( N\right) }$ starting from $\mathbf{\pi }_{0}$.

So with $X_{T_{\left( N\right) }}\overset{d}{\sim }\mathbf{\pi }_{\left(
N\right) }$, $X_{T_{\left( N\right) }}$ independent of $T_{\left( N\right) }$
and $\mathbf{P}\left( X_{n}=N\text{ for some }n<T_{\left( N\right) }\right)
=0$. The latter equation also indicates that $\tau _{\mathbf{\pi }%
_{0},N}\geq 1$ is stochastically larger than $\tau _{\mathbf{\pi }_{0},%
\mathbf{\pi }_{\left( N\right) }}\geq 1$ (statistically, $\left\{
X_{n}^{\left( N\right) }\right\} $ started from $\mathbf{\pi }_{0}$ enters $%
\mathbf{\pi }_{\left( N\right) }$ before first hitting state $N$)$.$

As a consequence,

\begin{proposition}
(\cite{Brown}, Corollary $4.1$) For all $n\geq 0$%
\begin{eqnarray*}
\text{sep}\left( \mathbf{P}_{\mathbf{\pi }_{0}}\left( X_{n}^{\left( N\right)
}=\cdot \right) ,\mathbf{\pi }_{\left( N\right) }\right)  &=&\underset{k}{%
\max }\left( 1-\mathbf{\pi }_{0}^{\prime }P_{\left( N\right) }^{n}\mathbf{e}%
_{k}/\pi _{\left( N\right) }\left( k\right) \right)  \\
&=&1-\mathbf{\pi }_{0}^{\prime }P_{\left( N\right) }^{n}\mathbf{e}_{N}/\pi
_{\left( N\right) }\left( N\right) =\mathbf{P}\left( T_{\left( N\right)
}>n\right) ,
\end{eqnarray*}
and $T_{\left( N\right) }$ is a minimal strong stationary time with
separating state $N$.
\end{proposition}

The separation distance sep$\left( \cdot ,\cdot \right) $ from $\mathbf{P}_{%
\mathbf{\pi }_{0}}\left( X_{n}^{\left( N\right) }=\cdot \right) $ to $%
\mathbf{\pi }_{\left( N\right) }$ gives an upper bound for the total
variation norm between these two probability measures.

\begin{equation*}
\mathbf{E}\left( T_{\left( N\right) }\right) =1+\sum_{n\geq 1}\frac{\mathbf{%
\pi }_{0}^{\prime }P_{\left( N\right) }^{n}\mathbf{e}_{N}}{\pi _{\left(
N\right) }\left( N\right) }=1+\frac{1}{\pi _{_{\left( N\right) }}\left(
N\right) }\mathbf{\pi }_{0}^{\prime }\left( I-P_{\left( N\right) }\right)
^{-1}P_{\left( N\right) }\mathbf{e}_{N}
\end{equation*}
There are some other facts pertaining to the fact that $\tau _{\mathbf{\pi }%
_{\left( N\right) },N}$ has a geometric convolution representation.

\begin{proposition}
\cite{Brown}

$\left( i\right) $ $\mathbf{P}\left( X_{n}^{\left( N\right) }=N\mid
X_{0}^{\left( N\right) }=N\right) $ is non-increasing with $n$, so 
\begin{equation}
\mathbf{P}\left( W_{1}^{\left( N\right) }>n\right) =\frac{P_{\left( N\right)
}^{n}\left( N,N\right) -\pi _{_{\left( N\right) }}\left( N\right) }{1-\pi
_{_{\left( N\right) }}\left( N\right) }.  \label{B40}
\end{equation}
is a well defined complementary mass function of some rv $W_{1}^{\left(
N\right) }$.

$\left( ii\right) $%
\begin{equation}
\tau _{\mathbf{\pi }_{\left( N\right) },N}=\sum_{i=1}^{G_{N}}W_{i}^{\left(
N\right) }  \label{B4}
\end{equation}
where $G_{N}\overset{d}{\sim }$geo$\left( \mathbf{\pi }_{_{\left( N\right)
}}\left( N\right) \right) $ (viz $\mathbf{P}\left( G_{N}=j\right) =\pi
_{_{\left( N\right) }}\left( N\right) \left( 1-\pi _{_{\left( N\right)
}}\left( N\right) \right) ^{j},$ $j=0,1,...$), independent of $W_{i}^{\left(
N\right) },$ $i\geq 1,$ an iid sequence with $W_{i}^{\left( N\right) }%
\overset{d}{=}W_{1}^{\left( N\right) }.$
\end{proposition}

\emph{Proof:} $\left( i\right) $ Because $P_{\left( N\right) }$ is SM as
well, $P_{\left( N\right) }^{n}\left( N,N\right) \geq P_{\left( N\right)
}^{n}\left( i,N\right) $ for all $i$ and $n$. Therefore, with $n_{2}>n_{1}$, 
\begin{equation*}
P_{\left( N\right) }^{n_{2}}\left( N,N\right) =\sum_{i=1}^{N}P_{\left(
N\right) }^{n_{2}-n_{1}}\left( N,i\right) P_{\left( N\right) }^{n_{1}}\left(
i,N\right) \leq P_{\left( N\right) }^{n_{1}}\left( N,N\right)
\sum_{i=1}^{N}P_{\left( N\right) }^{n_{2}-n_{1}}\left( N,i\right) =P_{\left(
N\right) }^{n_{1}}\left( N,N\right) .
\end{equation*}
As a result, $\mathbf{P}\left( X_{n}^{\left( N\right) }=N\mid X_{0}^{\left(
N\right) }=N\right) =\mathbf{e}_{N}^{\prime }P_{\left( N\right) }^{n}\mathbf{%
e}_{N}=P_{\left( N\right) }^{n}\left( N,N\right) $ is non-increasing with $n$
so that the law of $W_{1}^{\left( N\right) }$ is well-defined.

$\left( ii\right) $ 
\begin{equation*}
\sum_{n\geq 0}z^{n}\mathbf{P}\left( W_{1}^{\left( N\right) }>n\right) =\frac{%
1-\mathbf{E}\left( z^{W_{1}^{\left( N\right) }}\right) }{1-z}=\frac{1}{1-\pi
_{_{\left( N\right) }}\left( N\right) }\left( G_{N,N}\left( z\right) -\frac{%
\pi _{_{\left( N\right) }}\left( N\right) }{1-z}\right)
\end{equation*}
Using (\ref{B2}), we get 
\begin{equation}
\mathbf{E}\left( z^{\tau _{\mathbf{\pi }_{\left( N\right) },N}}\right) =%
\frac{1}{1+\frac{1-\pi _{_{\left( N\right) }}\left( N\right) }{\pi
_{_{\left( N\right) }}\left( N\right) }\left( 1-\mathbf{E}\left(
z^{W_{1}^{\left( N\right) }}\right) \right) }  \label{B5}
\end{equation}
which is the pgf of the geometric convolution $\sum_{i=1}^{G_{N}}W_{i}^{%
\left( N\right) }$. $\Box $

Note $G_{N}=0$ entails $\tau _{\mathbf{\pi }_{\left( N\right) },N}=0$, an
event with probability $\pi _{_{\left( N\right) }}\left( N\right) $.\newline

\emph{Remark:} Stochastically monotone Markov chain have a real and simple
second largest eigenvalue, \cite{Keilson}. Suppose $1=\lambda _{1}>\lambda
_{2}>\left| \lambda _{3}\right| \geq ...\geq \left| \lambda _{N}\right| >0$
where $\lambda _{k}=\lambda _{k,\left( N\right) }$ are the $N-$dependent
eigenvalues of $P_{\left( N\right) }$. Then, 
\begin{equation*}
\forall i,j\in \left\{ 1,...,N\right\} \text{, }\forall n\in \Bbb{N}\text{, }%
\exists c>0:\text{ }\left| P_{\left( N\right) }^{n}\left( i,j\right) -\pi
_{_{\left( N\right) }}\left( j\right) \right| \leq c\lambda _{2,\left(
N\right) }^{n}.
\end{equation*}
In particular, $\left| P_{\left( N\right) }^{n}\left( N,N\right) -\pi
_{_{\left( N\right) }}\left( N\right) \right| \leq c\lambda _{2,\left(
N\right) }^{n}$ and $P_{\left( N\right) }^{n}\left( N,N\right) $ is getting
close to $\pi _{_{\left( N\right) }}\left( N\right) $ as $n$ gets large,
useful for (\ref{B40}). $\diamondsuit $\newline

- \textbf{Quasi-stationary distribution (qsd).} An alternative point of view
on $\tau _{\mathbf{\pi }_{0},N}$ and $\tau _{\mathbf{\pi }_{\left( N\right)
},N}$ can also be seen from the classical theory of qsd's, \cite{cmsm}.

With $i\neq N$, let $\tau _{i,N}=\inf \left( n\geq 1:X_{n}^{\left( N\right)
}=N\mid X_{0}^{\left( N\right) }=i\right) $. We have 
\begin{equation*}
\mathbf{P}\left( \tau _{i,N}>1\right) =\mathbf{P}\left( X_{1}^{\left(
N\right) }\leq N-1\mid X_{0}^{\left( N\right) }=i\right) =P^{c}\left(
i,N-1\right) =F_{\left( N\right) }\left( N-1\right) ^{i}
\end{equation*}
\begin{eqnarray*}
\mathbf{P}\left( \tau _{i,N}>n+1\right) &=&\sum_{1\leq j<N}\mathbf{P}%
_{i}\left( X_{n}^{\left( N\right) }=j,\tau _{i,N}>n+1\right) \\
&=&\sum_{1\leq j<N}F_{\left( N\right) }\left( N-1\right) ^{j}\mathbf{P}%
_{i}\left( X_{n}^{\left( N\right) }=j,\tau _{i,N}>n\right)
\end{eqnarray*}
\begin{equation*}
\mathbf{P}\left( \tau _{i,N}>n+1\mid \tau _{i,N}>n\right) =\sum_{1\leq
j<N}F_{\left( N\right) }\left( N-1\right) ^{j}\mathbf{P}_{i}\left(
X_{n}^{\left( N\right) }=j\mid \tau _{i,N}>n\right)
\end{equation*}
\begin{equation*}
\underset{n\rightarrow \infty }{\rightarrow }\sum_{1\leq j<N}F_{\left(
N\right) }\left( N-1\right) ^{j}\mu _{\left( N-1\right) }\left( j\right) =:%
\mathbf{E}\left( z^{Z_{\left( N-1\right) }}\right) \mid _{z=F_{\left(
N\right) }\left( N-1\right) }=:\rho _{N}.
\end{equation*}
In the latter displayed formula, $\mathbf{\mu }_{\left( N-1\right) }\left(
\cdot \right) $ is the quasi-stationary limiting distribution of $\left\{
X_{n}^{\left( N\right) }\right\} $ when state $N$ has been removed and $%
Z_{\left( N-1\right) }\overset{d}{\sim }\mathbf{\mu }_{\left( N-1\right) }$.
Stated differently, $\mathbf{\mu }_{\left( N-1\right) }^{\prime }$ is the $%
\left( N-1\right) -$dimensional left eigenvector (associated to the dominant
eigenvalue $\rho _{N}<1$) of the substochastic matrix $P_{\left( N-1\right)
} $ obtained while removing the $N-$th row and column $N-$th column of $%
P_{\left( N\right) }$. We have used $\mathbf{P}_{i}\left( X_{n}^{\left(
N\right) }=j\mid \tau _{i,N}>n\right) \underset{n\rightarrow \infty }{%
\rightarrow }\mu _{\left( N-1\right) }\left( j\right) $, $j\in \left\{
1,...,N-1\right\} .$ Consequently,

\begin{proposition}
With $\rho _{N}$ the value of the pgf of $Z_{\left( N-1\right) }$ evaluated
at $F_{\left( N\right) }\left( N-1\right) $, independently of $i\in \left\{
1,...,N-1\right\} $%
\begin{equation*}
\underset{n\rightarrow \infty }{\lim }-\frac{1}{n}\log \mathbf{P}\left( \tau
_{i,N}>n\right) =-\log \mathbf{E}\left( z^{Z_{\left( N-1\right) }}\right)
\mid _{z=F_{\left( N\right) }\left( N-1\right) }=-\log \rho _{N}.
\end{equation*}
Equivalently, 
\begin{equation*}
\rho _{N}=\mathbf{E}\left( z^{Z_{\left( N-1\right) }}\right) \mid
_{z=F_{\left( N\right) }\left( N-1\right) }
\end{equation*}
is the rate of decay of $\mathbf{P}\left( \tau _{i,N}>n\right) $.
\end{proposition}

Similarly,

- With $\mathbf{\pi }_{0,0}$ defined by $\mathbf{\pi }_{0}^{\prime }=:\left( 
\mathbf{\pi }_{0,0}^{\prime },0\right) $, for any initial distribution $%
\mathbf{\pi }_{0,0},$%
\begin{equation*}
\underset{n\rightarrow \infty }{\lim }-\frac{1}{n}\log \mathbf{P}\left( \tau
_{\mathbf{\pi }_{0,0},N}>n\right) =-\log \mathbf{E}\left( z^{Z_{\left(
N-1\right) }}\right) \mid _{z=F_{\left( N\right) }\left( N-1\right) },
\end{equation*}
giving the decay rate of $\mathbf{P}\left( \tau _{\mathbf{\pi }%
_{0,0},N}>n\right) .$

- With $\mathbf{\pi }_{\left( N-1\right) }^{\prime }$ defined by $\mathbf{%
\pi }_{\left( N\right) }^{\prime }=\left( \mathbf{\pi }_{\left( N-1\right)
}^{\prime },\pi _{\left( N\right) }\left( N\right) \right) $, when starting
from the invariant measure 
\begin{equation*}
\underset{n\rightarrow \infty }{\lim }-\frac{1}{n}\log \mathbf{P}\left( \tau
_{\mathbf{\pi }_{\left( N-1\right) }^{\prime },N}>n\right) =-\log \mathbf{E}%
\left( z^{Z_{\left( N-1\right) }}\right) \mid _{z=F_{\left( N\right) }\left(
N-1\right) }.
\end{equation*}

- Clearly also, when the initial distribution coincides with the
quasi-stationary distribution: $\mathbf{\pi }_{0,0}=\mathbf{\mu }_{\left(
N-1\right) },$ 
\begin{equation*}
-\frac{1}{n}\log \mathbf{P}\left( \tau _{\mathbf{\mu }_{\left( N-1\right)
},N}>n\right) =-\log \mathbf{E}\left( z^{Z_{\left( N-1\right) }}\right) \mid
_{z=F_{\left( N\right) }\left( N-1\right) }=-\log \rho _{N}
\end{equation*}
for all $n.$ Letting 
\begin{equation*}
\mathbf{\mu }_{\left( N-1\right) }^{\prime }P_{\left( N-1\right) }=\rho _{N}%
\mathbf{\mu }_{\left( N-1\right) }^{^{\prime }}\text{ and }P_{\left(
N-1\right) }\mathbf{\phi }_{\left( N-1\right) }=\rho _{N}\mathbf{\phi }%
_{\left( N-1\right) },
\end{equation*}
be the $\left( N-1\right) $-dimensional left and right positive eigenvectors
of $P_{\left( N-1\right) }$ chosen so as to satisfy: $\left| \mathbf{\mu }%
_{\left( N-1\right) }\right| :=\sum_{j=1}^{N-1}\mu _{\left( N-1\right)
}\left( j\right) =1$ and $\mathbf{\mu }_{\left( N-1\right) }^{\prime }%
\mathbf{\phi }_{\left( N-1\right) }=1$, fixing the length $\left\| \mathbf{%
\phi }_{\left( N-1\right) }\right\| _{2}^{1/2}$ of $\mathbf{\phi }_{\left(
N-1\right) }$, then (by Perron-Frobenius theorem) 
\begin{equation*}
\rho _{N}^{-n}P_{\left( N-1\right) }^{n}\rightarrow \mathbf{\phi }_{\left(
N-1\right) }^{\prime }\mathbf{\mu }_{\left( N-1\right) }\text{ as }%
n\rightarrow \infty .
\end{equation*}
Hence, with $\mathbf{\pi }_{0}^{\prime }=\left( \mathbf{\pi }_{0,0}^{\prime
},0\right) ,$ $\left| \mathbf{\pi }_{0,0}\right| =1,$ and $\mathbf{\pi }%
_{\left( N\right) }^{\prime }=\left( \mathbf{\pi }_{\left( N-1\right)
}^{\prime },\pi _{\left( N\right) }\left( N\right) \right) $, $\left| 
\overline{\mathbf{\pi }}_{\left( N-1\right) }\right| <1$, and making use of $%
\tau _{N,N}=0$%
\begin{equation}
\mathbf{P}\left( \tau _{\mathbf{\pi }_{0},N}>n\right) =\mathbf{\pi }%
_{0,0}^{\prime }P_{\left( N-1\right) }^{n}\mathbf{1}\text{ and }\mathbf{P}%
\left( \tau _{\mathbf{\pi }_{\left( N\right) },N}>n\right) =\mathbf{\pi }%
_{\left( N-1\right) }^{\prime }P_{\left( N-1\right) }^{n}\mathbf{1<P}\left(
\tau _{\mathbf{\pi }_{0},N}>n\right)  \label{B6}
\end{equation}
\begin{equation*}
\mathbf{E}\left( \tau _{\mathbf{\pi }_{0},N}\right) =\mathbf{\pi }%
_{0,0}^{\prime }\left( I-P_{\left( N-1\right) }\right) ^{-1}\mathbf{1}\text{
and }\mathbf{E}\left( \tau _{\mathbf{\pi }_{\left( N\right) },N}\right) =%
\mathbf{\pi }_{\left( N-1\right) }^{\prime }\left( I-P_{\left( N-1\right)
}\right) ^{-1}\mathbf{1<E}\left( \tau _{\mathbf{\pi }_{0},N}\right)
\end{equation*}
and 
\begin{eqnarray*}
\rho _{N}^{-n}\mathbf{P}\left( \tau _{\mathbf{\pi }_{0},N}>n\right)
&\rightarrow &\mathbf{\pi }_{0,0}^{\prime }\mathbf{\phi }_{\left( N-1\right)
}\text{ as }n\rightarrow \infty \\
\rho _{N}^{-n}\mathbf{P}\left( \tau _{\mathbf{\pi }_{\left( N\right)
},N}>n\right) &\rightarrow &\mathbf{\pi }_{\left( N-1\right) }^{\prime }%
\mathbf{\phi }_{\left( N-1\right) }\text{ as }n\rightarrow \infty .
\end{eqnarray*}

\begin{proposition}
Suppose $\mathbf{\pi }_{0}$ is such that $\pi _{0}\left( i\right) /\pi
_{\left( N\right) }\left( i\right) $ decreases with $i$ and $\pi _{0}\left(
N\right) =0$. Then 
\begin{equation*}
\frac{\mathbf{P}\left( \tau _{\mathbf{\pi }_{0},N}>n\right) }{\mathbf{P}%
\left( \tau _{\mathbf{\pi }_{\left( N\right) },N}>n\right) }\underset{%
n\rightarrow \infty }{\rightarrow }\frac{\mathbf{\pi }_{0,0}^{\prime }%
\mathbf{\phi }_{\left( N-1\right) }}{\mathbf{\pi }_{\left( N-1\right)
}^{\prime }\mathbf{\phi }_{\left( N-1\right) }}\geq 1.
\end{equation*}
\end{proposition}

\emph{Proof:} Due to the stochastic domination of $\tau _{\mathbf{\pi }%
_{0},N}$ over $\tau _{\mathbf{\pi }_{\left( N\right) },N}$ stated in
Proposition $11$, the positive sequence 
\begin{equation*}
u_{n}:=\frac{\mathbf{P}\left( \tau _{\mathbf{\pi }_{0},N}>n\right) }{\mathbf{%
P}\left( \tau _{\mathbf{\pi }_{\left( N\right) },N}>n\right) }=\frac{\rho
_{N}^{-n}\mathbf{P}\left( \tau _{\mathbf{\pi }_{0},N}>n\right) }{\rho
_{N}^{-n}\mathbf{P}\left( \tau _{\mathbf{\pi }_{\left( N\right)
},N}>n\right) }
\end{equation*}
is bounded below by $1$ ($u_{n}\geq 1$ for all $n$). The sequence $u_{n}$ is
convergent with limit $u_{*}=\frac{\mathbf{\pi }_{0,0}^{\prime }\mathbf{\phi 
}_{\left( N-1\right) }}{\mathbf{\pi }_{\left( N-1\right) }^{\prime }\mathbf{%
\phi }_{\left( N-1\right) }}$ and the limit obeys $u_{*}\geq 1$.

We have $\rho _{N}^{-n}P_{\left( N-1\right) }^{n}\mathbf{1}\rightarrow 
\mathbf{\phi }_{\left( N-1\right) }$ as $n\rightarrow \infty $. The entries $%
\phi _{\left( N-1\right) }\left( i\right) $ are decreasing with $i$, because
it follows by induction that stochastic monotonicity of $P_{\left( N\right)
} $ implies the one of $P_{\left( N-1\right) }^{n}$, so that $\mathbf{e}%
_{i}^{\prime }P_{\left( N-1\right) }^{n}\mathbf{1}$ is decreasing with $i$.
Because $\pi _{0,0}\left( i\right) /\pi _{\left( N-1\right) }\left( i\right) 
$ is decreasing with $i$, the initial probability mass assigned to states
near the bottom state $\left\{ 1\right\} $ where $\mathbf{\phi }_{\left(
N-1\right) }$ takes its largest values exceeds the one assigned by $\mathbf{%
\pi }_{\left( N\right) }$. It is thus not that surprising that the numerator
of $u_{*}$ exceeds its denominator. $\Box $\newline

\emph{Remark:} From (\ref{B6}) 
\begin{eqnarray*}
\mathbf{E}\left( z^{\tau _{\mathbf{\pi }_{0},N}}\right) &=&1-\left(
1-z\right) \mathbf{\pi }_{0,0}^{\prime }\left( I-zP_{\left( N-1\right)
}\right) ^{-1}\mathbf{1}\text{ and } \\
\mathbf{E}\left( z^{\tau _{\mathbf{\pi }_{\left( N\right) },N}}\right)
&=&1-\left( 1-z\right) \mathbf{\pi }_{\left( N-1\right) }^{\prime }\left(
I-zP_{\left( N-1\right) }\right) ^{-1}\mathbf{1}
\end{eqnarray*}
Comparing the expression\emph{\ }of the pgf of\emph{\ }$\tau _{\mathbf{\pi }%
_{\left( N\right) },N}$ in terms of the Green kernel of $P_{\left(
N-1\right) }$ with (\ref{B5}), yields and identity for $\mathbf{\pi }%
_{\left( N-1\right) }^{\prime }\left( I-zP_{\left( N-1\right) }\right) ^{-1}%
\mathbf{1}$. Note that in (\ref{B5}), only the values of $\pi _{_{\left(
N\right) }}\left( N\right) $ and $P_{\left( N\right) }^{n}\left( N,N\right) $
matter. Comparing the expression\emph{\ }of the pgf of\emph{\ }$\tau _{%
\mathbf{\pi }_{0},N}$ in terms of the Green kernel of $P_{\left( N-1\right)
} $ with (\ref{B3a}), also yields and identity for $\mathbf{\pi }%
_{0,0}^{\prime }\left( I-zP_{\left( N-1\right) }\right) ^{-1}\mathbf{1}$. $%
\diamondsuit $\newline

By the definition of quasi-stationary distributions, we had 
\begin{equation*}
\mathbf{P}\left( X_{n}^{\left( N\right) }=j\mid \tau _{\mathbf{\pi }%
_{0},N}>n\right) \underset{n\rightarrow \infty }{\rightarrow }\mu _{\left(
N-1\right) }\left( j\right) \text{, }j\in \left\{ 1,...,N-1\right\} .
\end{equation*}
Because $P_{\left( N\right) }$ is stochastically monotone, Siegmund-Pollack
theorem holds, stating \cite{PS} 
\begin{equation*}
\mathbf{P}\left( X_{n}^{\left( N\right) }=j\mid \tau _{\mathbf{\pi }%
_{0},N}>n\right) \underset{n,N\rightarrow \infty }{\rightarrow }\pi \left(
j\right) \text{, }j\geq 1.
\end{equation*}
As $N$\ gets large, the qsd $\mathbf{\mu }_{\left( N-1\right) }$\ gets very
close to $\mathbf{\pi }_{\left( N-1\right) }$.\newline

- \textbf{Asymptotic exponentiality.}

- The rv $\tau _{\mathbf{\mu }_{\left( N-1\right) },N}$\ is geometric with
success parameter $1-\rho _{N}$, 
\begin{equation*}
\mathbf{E}\left( z^{\tau _{\mathbf{\mu }_{\left( N-1\right) },N}}\right) =%
\frac{z\left( 1-\rho _{N}\right) }{1-\rho _{N}z},
\end{equation*}
\ so with mean and variance $\mathbf{E}\left( \tau _{\mathbf{\mu }_{\left(
N-1\right) },N}\right) =1/\left( 1-\rho _{N}\right) $\ and $\sigma
^{2}\left( \tau _{\mathbf{\mu }_{\left( N-1\right) },N}\right) =\rho
_{N}/\left( 1-\rho _{N}\right) ^{2}.$ Suppose $\rho _{N}\rightarrow 1$ as $%
N\rightarrow \infty $. Then $\tau _{\mathbf{\pi }_{\left( N\right) },N}/%
\mathbf{E}\left( \tau _{\mathbf{\mu }_{\left( N-1\right) },N}\right) $\
becomes approximately exponential with mean $1$. We have $\mathbf{E}\left(
\tau _{\mathbf{\mu }_{\left( N-1\right) },N}\right) \rightarrow \infty $
while $\sigma \left( \tau _{\mathbf{\mu }_{\left( N-1\right) },N}\right) /%
\mathbf{E}\left( \tau _{\mathbf{\mu }_{\left( N-1\right) },N}\right) =\sqrt{%
\rho _{N}}\rightarrow 1$, as\emph{\ }$N\rightarrow \infty $\emph{.}

- Brown raised the question of asymptotic exponentiality of $\tau _{\mathbf{%
\pi }_{\left( N\right) },N}/\mathbf{E}\left( \tau _{\mathbf{\pi }_{\left(
N\right) },N}\right) $.

If $\sigma ^{2}\left( \tau _{\mathbf{\pi }_{\left( N\right) },N}\right)
<\infty $, as a scaled geometric convolution, $\tau _{\mathbf{\pi }_{\left(
N\right) },N}/\mathbf{E}\left( \tau _{\mathbf{\pi }_{\left( N\right)
},N}\right) $ is approximately exponential if $\mathbf{E}\left( \tau _{%
\mathbf{\pi }_{\left( N\right) },N}\right) \rightarrow \infty $ while $%
\sigma \left( \tau _{\mathbf{\pi }_{\left( N\right) },N}\right) /\mathbf{E}%
\left( \tau _{\mathbf{\pi }_{\left( N\right) },N}\right) \rightarrow 1$, as%
\emph{\ }$N\rightarrow \infty $ for the truncated Lamperti model with
truncated target distribution $\mathbf{\pi }_{\left( N\right) }$. Error
bounds can be obtained from the first two moments of $\tau _{\mathbf{\pi }%
_{\left( N\right) },N}$ given by

\begin{eqnarray*}
\mathbf{E}\left( \tau _{\mathbf{\pi }_{\left( N\right) },N}\right) &=&%
\mathbf{E}\left( G_{N}\right) \mathbf{E}\left( W_{1}^{\left( N\right)
}\right) =\frac{1-\pi _{_{\left( N\right) }}\left( N\right) }{\pi _{_{\left(
N\right) }}\left( N\right) }\mathbf{E}\left( W_{1}^{\left( N\right) }\right)
\\
&=&\frac{1}{\pi _{_{\left( N\right) }}\left( N\right) }\sum_{n\geq 0}\left(
P_{\left( N\right) }^{n}\left( N,N\right) -\pi _{_{\left( N\right) }}\left(
N\right) \right)
\end{eqnarray*}
\begin{eqnarray*}
\sigma ^{2}\left( \tau _{\mathbf{\pi }_{\left( N\right) },N}\right) &=&%
\mathbf{E}\left( G_{N}\right) \sigma ^{2}\left( W_{1}^{\left( N\right)
}\right) +\left( \mathbf{E}\left( W_{1}^{\left( N\right) }\right) \right)
^{2}\sigma ^{2}\left( G,_{N}\right) \\
\mathbf{E}\left( \tau _{\mathbf{\pi }_{\left( N\right) },N}^{2}\right) &=&2%
\mathbf{E}\left( \tau _{\mathbf{\pi }_{\left( N\right) },N}\right) ^{2}+%
\mathbf{E}\left( G_{N}\right) \mathbf{E}\left( \left( W_{1}^{\left( N\right)
}\right) ^{2}\right) .
\end{eqnarray*}
The question of the approximation by an exponential distribution also arises
for $\tau _{\mathbf{\pi }_{0},N}/\mathbf{E}\left( \tau _{\mathbf{\pi }%
_{0},N}\right) $. In this direction indeed,

\begin{proposition}
(\cite{Brown}, \cite{Brown2}) With $t\geq 0$
\end{proposition}

\begin{eqnarray*}
\sup_{t}\left| \mathbf{P}\left( \frac{\tau _{\mathbf{\pi }_{\left( N\right)
},N}}{\mathbf{E}\left( \tau _{\mathbf{\pi }_{\left( N\right) },N}\right) }%
>t\right) -e^{-t}\right| &\leq &\pi _{_{\left( N\right) }}\left( N\right) 
\frac{\mathbf{E}\left( \left( W_{1}^{\left( N\right) }\right) ^{2}\right) }{%
\mathbf{E}\left( W_{1}^{\left( N\right) }\right) ^{2}} \\
&=&2\left( 1-\pi _{_{\left( N\right) }}\left( N\right) \right) \left[ \frac{%
\mathbf{E}\left( \tau _{\mathbf{\pi }_{\left( N\right) },N}^{2}\right) }{2%
\mathbf{E}\left( \tau _{\mathbf{\pi }_{\left( N\right) },N}\right) ^{2}}%
-1\right] , \\
\sup_{t}\left| \mathbf{P}\left( \frac{\tau _{\mathbf{\pi }_{0},N}}{\mathbf{E}%
\left( \tau _{\mathbf{\pi }_{0},N}\right) }>t\right) -e^{-t}\right| &\leq &%
\frac{\mathbf{E}\left( T_{\left( N\right) }\right) }{\mathbf{E}\left( \tau _{%
\mathbf{\pi }_{\left( N\right) },N}\right) }+2\left( 1-\pi _{_{\left(
N\right) }}\left( N\right) \right) \left[ \frac{\mathbf{E}\left( \tau _{%
\mathbf{\pi }_{\left( N\right) },N}^{2}\right) }{2\mathbf{E}\left( \tau _{%
\mathbf{\pi }_{\left( N\right) },N}\right) ^{2}}-1\right]
\end{eqnarray*}
gives the sup-norm distance between respectively $\tau _{\mathbf{\pi }%
_{\left( N\right) },N}/\mathbf{E}\left( \tau _{\mathbf{\pi }_{\left(
N\right) },N}\right) $, $\tau _{\mathbf{\pi }_{0},N}/\mathbf{E}\left( \tau _{%
\mathbf{\pi }_{0},N}\right) $ and an exponential rv with mean $1.$

This shows that if, as $N$ grows large, the mean and standard deviation of $%
\tau _{\mathbf{\pi }_{\left( N\right) },N}/\mathbf{E}\left( \tau _{\mathbf{%
\pi }_{\left( N\right) },N}\right) $ behave like the one of an exponential
distribution that is if $\sigma \left( \tau _{\mathbf{\pi }_{\left( N\right)
},N}\right) /\mathbf{E}\left( \tau _{\mathbf{\pi }_{\left( N\right)
},N}\right) \rightarrow 1$, then $\mathbf{E}\left( \tau _{\mathbf{\pi }%
_{\left( N\right) },N}^{2}\right) /\left( 2\mathbf{E}\left( \tau _{\mathbf{%
\pi }_{\left( N\right) },N}\right) ^{2}\right) \rightarrow 1$ and the
exponential approximation for the law of $\tau _{\mathbf{\pi }_{\left(
N\right) },N}/\mathbf{E}\left( \tau _{\mathbf{\pi }_{\left( N\right)
},N}\right) $ is valid. If in addition, as $N$ becomes large 
\begin{equation*}
\mathbf{E}\left( T_{\left( N\right) }\right) /\mathbf{E}\left( \tau _{%
\mathbf{\pi }_{\left( N\right) },N}\right) \ll 2\left( 1-\pi _{_{\left(
N\right) }}\left( N\right) \right) \left[ \frac{\mathbf{E}\left( \tau _{%
\mathbf{\pi }_{\left( N\right) },N}^{2}\right) }{2\mathbf{E}\left( \tau _{%
\mathbf{\pi }_{\left( N\right) },N}\right) ^{2}}-1\right] ,
\end{equation*}
then the same holds true for the law of $\tau _{\mathbf{\pi }_{0},N}/\mathbf{%
E}\left( \tau _{\mathbf{\pi }_{0},N}\right) $.\newline

- \textbf{Time reversal. }Consider the time-reversed\textbf{\ }version $%
X_{n,\left( N\right) }^{\leftarrow }$ of the truncated Lamperti chain, so
with one-step transition matrix 
\begin{equation*}
\overleftarrow{P}_{\left( N\right) }=D_{\mathbf{\pi }_{\left( N\right)
}}^{-1}P_{\left( N\right) }^{\prime }D_{\mathbf{\pi }_{\left( N\right) }}.
\end{equation*}
Its time-reversed transition matrix being $P_{\left( N\right) }$ which is in
particular stochastically monotone, the Brown theory for hitting times
applies to the time-reversed process as well (see \cite{Brown}), with $%
\overleftarrow{\tau }_{\mathbf{\pi }_{0},N}$ and $\overleftarrow{\tau }_{%
\mathbf{\pi }_{\left( N\right) },N}$ standing for the hitting times of the
time-reversed chain. The time-reversed process $X_{n,\left( N\right)
}^{\leftarrow }$ thus constructed is a truncated version of the process
defined from (\ref{TR}).\newline

\textbf{Acknowledgments:} The authors are indebted for support of the
CMM-Basal Conicyt project AFB170001 and I.E.A. Cergy. T. H. acknowledges
partial support from the labex MME-DII (Mod\`{e}les Math\'{e}matiques et
\'{E}conomiques de la Dynamique, de l' Incertitude et des Interactions),
ANR11-LBX-0023-01. This work also benefited from the support of the Chair
``Mod\'{e}lisation Math\'{e}matique et Biodiversit\'{e}'' of Veolia-Ecole
Polytechnique-MNHN-Fondation X.

\end{document}